\documentclass[11pt]{article}
\usepackage{ulem}
\usepackage{graphicx,amsmath,amsfonts,amssymb,color,bbm}
\usepackage{epsfig}

  \newcommand{\splt}[1]{{\mbox{\rm split} #1}}

 \newcommand{\pend}{\hfill \thicklines \framebox(5.5,5.5)[l]{}}
 \newenvironment{proof}{\noindent {\sc  Proof.} \rm}{\pend}

\numberwithin{equation}{section}

 \newtheorem{theorem}{Theorem}[section]
 \newtheorem{lemma}{Lemma}[section]
 
 \newtheorem{proposition}{Proposition}[section]
 
 \newtheorem{remark}{Remark}[section]

 \newtheorem{definition}{Definition}[section]

\begin{document}
 \pagenumbering{arabic} \thispagestyle{empty}
\setcounter{page}{1}

\title
{Tail Asymptotics for the $M_1,M_2/G_1,G_2/1$ Retrial Queue with Priority}
\author{
    Bin Liu $^{1, a}$ and Yiqiang Q. Zhao $^b$
 \\ {\small a. School of Mathematics and Physics, Anhui Jianzhu University,
Hefei 230601, P.R. China}\\
{\small b. School of Mathematics and Statistics,
Carleton University, Ottawa, ON, Canada K1S 5B6}}
\date{\today}

\footnotetext[1]{Corresponding author: Bin Liu, E-mail address:
bliu@amt.ac.cn}
\maketitle

\begin{abstract}
	
Stochastic networks with complex structures are key modelling tools for many important applications.
In this paper, we consider a specific type of network: the retrial queueing systems with priority. This type of queueing system is important in various applications, including telecommunication and computer management networks with big data. For this type of system, we propose a detailed stochastic decomposition approach to
study its asymptotic behaviour of the tail probability of the number of customers in the steady-state for retrial queues
with two types (Type-1 and Type-2) of customers, in which Type-1 customers (in a queue) have non-preemptive priority to receive service over  Type-2 customers (in an orbit). Under the assumption that the service times of Type-1 customers have a regularly varying tail and the service times of Type-2 customers have a tail lighter than Type-1 customers, we obtain tail asymptotic properties for the number of customers in the queue and in the orbit, respectively, conditional on the server's status, in terms of a detailed stochastic decomposition approach. Tail asymptotic properties are often used as key tools for approximating various performance metrics and constructing numerical algorithms for computations.

\medskip

\medskip

\end{abstract}

\section{Introduction} \label{sec:1}

Rapid advances in the fields of computer and communication technologies, with fast increasing internet, big data and smart phone applications, have significantly changed every aspect of our life. These accelerated developments have continuously raised new challenges in modelling, performance analysis, system control and optimization. As a consequence of these challenges, the resulting stochastic networks, as key modelling tools, become progressively complex, due to dependence structures, dimensions, and the size of the data involved. For such networks, exact solutions are often rare, whereas asymptotic behaviours and properties are among the key candidates that we search for.
We consider a single-server retrial queue with two types of customers (Type-1 and Type-2), denoted by $M_1,M_2/G_1,G_2/1$. This model was studied by Falin, Artalejo and Martin in \cite{Falin-Artalejo-Martin:1993}. In this model,
customers arrive according to a Poisson process at rate $\lambda>0$ and with probabilities $q\in (0,1)$ and $p=1-q$ to be Type-1 and Type-2, respectively. In other words, Type-1 and Type-2 customers form two independent Poisson arrival processes with rates $\lambda_1\equiv\lambda q$ and $\lambda_2\equiv\lambda p$, respectively. If the server is idle upon the arrival of a Type-1 or Type-2 customer, the customer receives the service immediately and leaves the system after the completion of service. If an arriving Type-1 customer finds the server being busy, it joins the priority queue with an infinite waiting capacity. If a Type-2 customer finds the server being busy upon arrival, it enters the orbit and make retrial attempts later for receiving a service. Each of the Type-2 customers in the orbit repeatedly tries, independently, to receive service according to a Poisson process with a common retrial rate $\mu$ until it finds the server being idle, and receives its service immediately. Type-1 customers have non-preemptive priority to receive service over  Type-2 customers. Thus, as long as the priority queue is not empty, all retrials by Type-2 customers from the orbit are blocked (or failed), and all
blocked Type-2 customers return to the orbit with probability one. Type-$i$ customers have service time $T_{\beta_i}$, whose probability distribution is $F_{\beta_i}(x)$ with $F_{\beta_i}(0)=0$, and $T_{\beta_i}$ is assumed to have a finite mean $\beta_{i,1}$, $i=1,2$, where the second subscript is used to indicate the first moment of the service time. The Laplace-Stieltjes transforms (LST) of distribution function $F_{\beta_i}(x)$ is denoted by $\beta_i(s)$, $i=1,2$.
Let $\rho_1=\lambda_1\beta_{1,1}$, $\rho_2=\lambda_2\beta_{2,1}$ and $\rho=\rho_1+\rho_2=\lambda(q\beta_{1,1}+p\beta_{2,1})$. It follows from \cite{Falin-Artalejo-Martin:1993} that this system is stable if and only if $\rho<1$. We will assume that $\rho<1$ throughout this paper.

We refer readers to the following books, or review articles, for an updated status of studies on retrial queues and for more references therein:
Falin~\cite{Falin:1990},  Artalejo and G\'{o}mez-CorralFalin~\cite{Artalejo-GomezCorral:2008}, Kim and Kim~\cite{Kim-Kim:2016}, and
Phung-Duc~\cite{Phung-Duc:2017}. We also mention here the following two references, which are closely related to the study in this paper:
Kim, Kim and Ko~\cite{Kim-Kim-Ko:2007}, and Kim, Kim and Kim~\cite{Kim-Kim-Kim:2010c}.
Priority retrial queueing systems are a type of very important retrial queues, which find many applications, for example, in computer network management and telecommunication systems. In such systems, there are usually two or more types of customers. A survey of studies on single server retrial queues with priority calls (or customers), published by 1999, can be found in Choi and Chang~\cite{Choi-Chang:1999}. Since then, more publications on priority retrial queues are available, such as
Artalejo, Dudin and Klimenok~\cite{Artalejo-Dudin-Klimenok:2001},
Lee~\cite{Lee:2001},
G\'{o}mea-Corral~\cite{Gomez-Corral:2002},
Wang~\cite{Wang:2008},
Dimitriou~\cite{Dimitriou-2013b},
Wu and Lian~\cite{Wu-Lian:2013},
Wu, Wang and Liu~\cite{Wu-Wang-Liu:2013},
Gao~\cite{Gao:2015},
Dudin \textit{et al.}~\cite{Dudin-etc:2017},
Walraevens, Claeys and Phung-Duc~\cite{Walraevens-Claeys-Phung-Duc:2018}, among possible others.
Readers may refer to \cite{Dimitriou-2013b,Wu-Wang-Liu:2013} for more detailed reviews of the above mentioned studies.

Different from the above mentioned studies, our focus in this paper is on heavy-tailed behaviour of stationary (conditional) probabilities (assuming the stability of the system). Specifically, we assume that the tail probability of the service time for Type-1 customers is regularly varying, and the tail probability of the service time for Type-2 customers is lighter than that for Type-1 customers (see Assumptions~A1 and A2). Under these assumptions, we characterize the tail asymptotic behaviour for the following key system performance metrics:
\begin{description}
\item[PO-0] Conditional tail probability of the number of customers in the orbit given that the server is idle;
\item[PO-1] Conditional tail probability of the number of customers in the orbit given that the server is serving a Type-1 customer;
\item[PO-2] Conditional tail probability of the number of customers in the orbit given that the server is serving a Type-2 customer;
\item[PQ-1] Conditional tail probability of the number of customers in the queue given that the server is serving a Type-1 customer;
\item[PQ-2] Conditional tail probability of the number of customers in the queue given that the server is serving a Type-2 customer.
\end{description}
It is obvious that the queue should be empty when the server is idle. The tail asymptotic behaviour is one of the key subjects in applied probability. It is also very useful in approximations and computations, such as providing performance metrics and developing numerical algorithms (see Liu and Zhao~\cite{Liu-Zhao:2010} for some of its applications).

The main discovery in this paper is that the tail for all of the above mentioned conditional probabilities is also regularly varying with a dominant influence by the service time distribution for Type-1 customers, except for PQ-2, the tail of which is dominated by the service time for Type-2 customers (see Theorems~\ref{the:PO-0} and \ref{The:4.2} for details).
To obtain our main result, we propose a detailed stochastic decomposition approach, which has been recently applied for tail asymptotic analysis in various queueing models, including Liu, Wang and Zhao~\cite{Liu-Wang-Zhao:2014, Liu-Wang-Zhao:2012}, Liu, Min and Zhao~\cite{Liu-Min-Zhao:2017}, and  Liu and Zhao~\cite{Liu-Zhao:2017b, Liu-Zhao:2018}. Stochastic decomposition has been widely used in queueing system analysis. For example, it is well known that for the $M/G/1$ retrial queue, one can stochastically decompose the total number of customers in the system as the independent sum of the total number of customers in the corresponding standard (without retrials) $M/G/1$ queueing system and another random variable. The detailed stochastic decomposition approach is also to decompose a random variable, for example the number of customers in the queue, into a sum of independent variables, but with more detail. In the detailed decomposition, the sum consists of a fixed, or random, number of independent random variables (summands) such that the tail asymptotic property for each  summand is available, and a detailed analysis allows us to identify the summands, which play a dominant role for the tail asymptotic behaviour of the random sum.

The rest of the paper is organized as follows: In Section~\ref{sec:2}, we provide expressions for the probability generating functions of interest, which are our starting point. In Section~\ref{sec:3}, detailed stochastic decompositions are obtained. In Section~\ref{sec:4}, tail asymptotic properties, for each of the decomposed components, are discussed, which lead to our main results. This section also contains a concluding remark. Most of the literature results, needed in this paper, are collected in the appendix.

\section{Preliminary} \label{sec:2}

Let $R_{que}$ be the number of Type-1 customers in the queue, {\it excluding} the possible one in the service, let $R_{orb}$ be the number of Type-2 customers in the orbit, and let $I_{ser}=0,1\mbox{ or }2$ according to the status of the server: idle,  busy with a Type-1 customer, or busy with a Type-2 customer, respectively. Let $R_0$ be a random variable (r.v.), whose distribution coincides with the conditional distribution of $R_{orb}$ given that $I_{ser}=0$, and let $(R_{11},R_{12})$ and $(R_{21},R_{22})$ be  two-dimensional r.v.s, whose distributions coincide with the conditional distributions of $(R_{que},R_{orb})$ given that $I_{ser}=1$ and $I_{ser}=2$, respectively. Precisely,  $R_0$ has the probability generating function (PGF): $R_0(z_2)=E(z_2^{R_0})\overset{\footnotesize\mbox{def}}{=}E(z_2^{R_{orb}}|I_{ser}=0)$, and $(R_{i1},R_{i2})$ has the PGF: $R_i(z_1,z_2)=E(z_1^{R_{i1}}z_2^{R_{i2}})\overset{\footnotesize\mbox{def}}{=}E(z_1^{R_{que}}z_2^{R_{orb}}|I_{ser}=i)$ for $i=1,2$.

The following expressions for $R_0(z_2)$, $R_1(z_1,z_2)$ and $R_2(z_1,z_2)$ were obtained by Falin, Artalejo and Martin~\cite{Falin-Artalejo-Martin:1993}, which will be our starting point for tail asymptotics:
$P\{I_{ser}=0\}=1-\rho$, $P\{I_{ser}=1\}=\rho_1$, $P\{I_{ser}=2\}=\rho_2$,
\begin{align}
R_0(z_2) =& \exp\left\{-\frac {\lambda} {\mu}\int_{z_2}^{1}
\frac {1 - \beta(\lambda-\lambda g(u))} {\beta_2(\lambda-\lambda g(u))-u} du\right\},\label{P0-1}\\
R_1(z_1,z_2) =& \frac {1-\rho} {\rho_1}\cdot \frac {W(z_1,z_2)} {(\beta_2(\lambda-\lambda g(z_2))-z_2) (z_1-\beta_1(\lambda-\lambda_1 z_1-\lambda_2 z_2))}\nonumber\\
 &\cdot \frac {1-\beta_1(\lambda-\lambda_1 z_1 -\lambda_2 z_2)} {\lambda-\lambda_1 z_1 -\lambda_2 z_2}\cdot R_0(z_2),\label{P1-1}
\end{align}
and
\begin{equation}
R_2(z_1,z_2)=\frac {1-\rho} {\rho_2}\cdot\frac {\lambda-\lambda g(z_2)} {\beta_2(\lambda-\lambda g(z_2))-z_2}\cdot \frac {1-\beta_2(\lambda-\lambda_1 z_1-\lambda_2 z_2)} {\lambda-\lambda_1 z_1-\lambda_2 z_2}\cdot R_0(z_2), \label{P2-1}
\end{equation}
where
\begin{align}
\beta(s) =&q\beta_1(s) + p\beta_2(s),\label{beta(s)-1}\\
g(z_2) =&q h(z_2) + p z_2 \label{g(z_2)-1}
\end{align}
and
\begin{align}
W(z_1,z_2) =&[\lambda-\lambda_1 h(z_2)-\lambda_2 \beta_2(\lambda-\lambda g(z_2))] (\beta_2(\lambda-\lambda_1 z_1-\lambda_2 z_2)-z_2)\nonumber\\
& +(\lambda-\lambda_1 z_1 -\lambda_2 \beta_2(\lambda-\lambda_1 z_1-\lambda_2 z_2)) [z_2-\beta_2(\lambda-\lambda g(z_2))] \label{W-1}
\end{align}
with the function $h(z_2)$ being determined uniquely by the equation
\begin{equation}
h(z_2)=\beta_1(\lambda-\lambda_1 h(z_2)-\lambda_2 z_2). \label{Falin-phi-eqn}
\end{equation}

While we are not expecting to have any closed formulas for the inverse functions (or probabilities) of the above transformation functions, it is our focus in this paper to use the stochastic decomposition ideas to obtain simple characterizations of the tail probabilities. This technique is referred to as the {\it detailed stochastic decomposition approach} for transformation functions. To this end,
it is worth mentioning that (i) $\beta(s)$ in (\ref{beta(s)-1}) is the LST of the mixed distribution $F_{\beta}(x)\stackrel{\rm def}{=}q F_{\beta_1}(x)+p F_{\beta_2}(x)$; (ii) both $h(z_2)$ in (\ref{Falin-phi-eqn}) and $g(z_2)$  in (\ref{g(z_2)-1}) can be regarded as the PGFs of r.v.s, which will be verified in the next subsection.

\subsection{Probabilistic interpretations for PGF $h(z_2)$ and $g(z_2)$}

We will show that $h(z_2)$ is closely related to the busy period $T_{\alpha}$ of the standard $M/G/1$ queue with arrival rate $\lambda_1$ and the service time $T_{\beta_1}$.
By $F_{\alpha}(x)$ we denote the probability distribution function of $T_{\alpha}$, and by $\alpha(s)$ the LST of $F_{\alpha}(x)$.
The following are well-known results about this $M/G/1$ queue:
\begin{eqnarray}
\alpha(s) &=&\beta_1(s+\lambda_1- \lambda_1 \alpha(s)),\label{busy-eqn-alpha}\\
\alpha_1 &\stackrel{\rm def}{=}&E(T_{\alpha})=\beta_{1,1}/(1-\rho_1).\label{busy-exp-T-alpha}
\end{eqnarray}

Throughout this paper, we will use the notation $N_b(t)$ to represent the number of Poisson arrivals with rate $b$ within the time interval $(0,t]$.
Now, let us consider $N_{\lambda_2}(T_{\alpha})$, the number of arrivals of a Poisson process at rate $\lambda_2$ within an independent random time $T_{\alpha}$.
The PGF of $N_{\lambda_2}(T_{\alpha})$ is easily obtained as follows:
\begin{equation}\label{alpha-phi-0}
E(z_2^{N_{\lambda_2}(T_{\alpha})})=\int_0^{\infty}\sum_{n=0}^{\infty}z_2^n\frac {(\lambda_2 x)^n} {n!} e^{-\lambda_2 x}dF_{\alpha}(x)=\alpha(\lambda_2-\lambda_2 z_2).
\end{equation}
It follows from (\ref{busy-eqn-alpha}) that
\begin{equation}
\alpha(\lambda_2-\lambda_2 z_2)=\beta_1(\lambda - \lambda_1 \alpha(\lambda_2-\lambda_2 z_2) -\lambda_2 z_2).\label{busy-eqn-alpha-z}
\end{equation}
By comparing (\ref{Falin-phi-eqn}) and (\ref{busy-eqn-alpha-z}) and noticing the uniqueness of $h(z_2)$, we immediately have
\begin{equation}
h(z_2)=\alpha(\lambda_2-\lambda_2 z_2),\label{alpha-phi}
\end{equation}
which, together with (\ref{alpha-phi-0}), implies that $h(z_2)=E(z_2^{N_{\lambda_2}(T_{\alpha})})$ is the PGF of the non-negative integer-valued r.v. $N_{\lambda_2}(T_{\alpha})$, which is the number of Poisson arrivals, with arrival rate $\lambda_2$, during a busy period for the standard $M/G/1$ queue with arrival rate $\lambda_1$ and service time $T_{\beta_1}$. In addition, $g(z_2)$, as defined in (\ref{g(z_2)-1}), is also a PGF of non-negative integer-valued r.v., denoted by $X_g$,
i.e., $g(z_2)=E(z_2^{X_g})$.
It follows from (\ref{g(z_2)-1}) that
\begin{equation}\label{G-def}
X_g\stackrel{\rm d}{=}\left\{\begin{array}{ll}
1 &\mbox{ with probability }p,\\
N_{\lambda_2}(T_{\alpha}), &\mbox{ with probability }q,
\end{array}
\right.
\end{equation}
where we have used the symbol
$``\stackrel{\rm d}{=}"$ to mean the equality in probability distribution. Such a symbol will be used throughout the paper.

It is easy to obtain that $E(N_{\lambda_2}(T_{\alpha}))=\lambda_2\beta_{1,1}/(1-\rho_1)$ and
\begin{equation}\label{G-exp}
E(X_g)=p+q\lambda_2\beta_{1,1}/(1-\rho_1)=p/(1-\rho_1).
\end{equation}

\subsection{Assumptions on service times}

%

It is well known that for a distribution $F$ on $(0,\infty)$, if $\overline{F}$ is regularly varying with index $-\sigma$, $\sigma\ge 0$ (see Definition~\ref{Definition 3.1}) or $\overline{F}\in \mathcal{R}_{-\sigma}$, then $F$ is subexponnetial (see Definition~\ref{Definition 3.2}) or $F\in \mathcal S$
(see, e.g., Embrechts et al. \cite{Embrechts1997}).
We will use $L(t)$ to represent a slowly varying function at $\infty$ and make the following basic assumptions on the service time $T_{\beta_i}$ of Type-$i$ customers, $i=1,2$.
\begin{description}
\item[A1.] {\it $T_{\beta_1}$ has tail probability
$P\{T_{\beta_1}>t\} \sim t^{-a_1}L(t)$ as $t\to\infty$, where $a_1>1$.}
\item[A2.] {\it $T_{\beta_2}$ has tail probability
$P\{T_{\beta_2}>t\} \sim e^{-r t} t^{-a_2} L(t)$ as $t\to\infty$, where $-\infty<a_2<\infty$ if $r > 0$, or $a_2>a_1$ if $r=0$.}
\end{description}

Clearly, under assumptions A1 and A2, the service time $T_{\beta_1}$ of Type-1 customers has a tail probability heavier than the service time $T_{\beta_2}$ of Type-2 customers. If $r>0$, $T_{\beta_2}$ has a light tail, i.e., $E(e^{\varepsilon T_{\beta_2}})<\infty$ for some $\varepsilon>0$. If $r=0$, then $T_{\beta_2}$ has a regularly varying tail with index $-a_2$.

Since $T_{\alpha}$ is the busy period of the ordinary $M/G/1$ queue with arrival rate $\lambda_1$ and the service time $T_{\beta_1}$, its asymptotic tail probability is regularly varying according to de Meyer and Teugels~\cite{Meyer-Teugels:1980} (see Lemma~\ref{Lemma 2.1} in Appendix).)

\section{Detailed stochastic decompositions} \label{sec:3}

In this section, we will apply the detailed stochastic decomposition technique to
r.v.s $R_0$, $(R_{11},R_{12})$ and $(R_{21},R_{22})$. The decomposition results obtained will be used in asymptotic analysis later in Section~\ref{sec:4}.
First, we rewrite (\ref{P0-1}). Let
\begin{eqnarray}
K_a(u)&=&\frac {1-\rho_1} {p}\cdot\frac {1- g(u)} {1-u},\label{Ka-1}\\
K_b(u)&=&\frac 1 {\rho}\cdot\frac {1- \beta (\lambda- \lambda g(u))} {1- g(u)},\label{Kb-1}\\
K_c(u)&=&\frac{1-\rho} {1-\rho_1}\cdot\frac {1-u} {\beta_2(\lambda-\lambda g(u))-u}.\label{Kc-1}
\end{eqnarray}
Immediately, we have,
\begin{eqnarray}\label{K=Kabc}
\frac {1 - \beta(\lambda-\lambda g(u))} {\beta_2(\lambda-\lambda g(u))-u}
&=&\frac {\rho p} {1-\rho}\cdot K_a(u)\cdot K_b(u)\cdot K_c(u).
\end{eqnarray}
Substituting (\ref{K=Kabc}) into (\ref{P0-1}), we obtain
\begin{equation}\label{D0-2}
R_0(z_2)=\exp\left\{-\psi\int_{z_2}^{1}
K(u) du\right\},
\end{equation}
where
\begin{eqnarray}
\psi&=&\rho\lambda_2/(\mu (1-\rho)),\label{psi}\\
K(u)&=&K_a(u)\cdot K_b(u)\cdot K_c(u).\label{K(u)-1}
\end{eqnarray}
In the next subsection, we will verify that $K_a(u)$, $K_b(u)$, $K_c(u)$ and $K(u)$ can be viewed as the PGFs of four r.v.s, denoted by $K_a$, $K_b$, $K_c$ and $K$, respectively.

Let $F_{\alpha}^{(e)}(x)$ be the so-called equilibrium distribution of $F_{\alpha}(x)$,
which is defined as $F_{\alpha}^{(e)}(x)= \alpha_1^{-1}\int_0^{x}(1-F_{\alpha}(t))dt$, where $\alpha_1=E(T_{\alpha})$ given in (\ref{busy-exp-T-alpha}).
The LST of $F_{\alpha}^{(e)}(x)$ can be written as $\alpha^{(e)}(s)=(1-\alpha(s))/(\alpha_1 s)$. Similarly,
$F_{\beta_i}^{(e)}(x) \stackrel{\rm def}{=} \beta_{i,1}^{-1}\int_0^{x}(1-F_{\beta_i}(t))dt$, $F_{\beta}^{(e)}(x) \stackrel{\rm def}{=} (q\beta_{1,1}+p\beta_{2,1})^{-1}\int_0^{x}(1-F_{\beta}(t))dt$, and the LSTs of $F_{\beta_i}^{(e)}(x)$ and $F_{\beta}^{(e)}(x)$ can be written as $\beta_i^{(e)}(s)=(1-\beta_i(s))/(\beta_{i,1} s)$, $i=1,2$, and $\beta^{(e)}(s)=(1-\beta(s))/((q\beta_{1,1}+p\beta_{2,1}) s)$, respectively.

\subsection{Stochastic decomposition on $K$}

By (\ref{alpha-phi}) and the definition of $\alpha^{(e)}(s)$, we can write $(1- h(u))/(1-u)= \lambda_2\alpha_1\cdot \alpha^{(e)}(\lambda_2 - \lambda_2 u)$, from which, and by (\ref{Ka-1}), (\ref{g(z_2)-1}) and (\ref{busy-exp-T-alpha}), we have,
\begin{eqnarray}\label{Ka(u)}
K_a(u)&=&\frac {1-\rho_1} {p}\left [q\cdot\frac {1-h(u)} {1-u} + p\right]=\rho_1 \alpha^{(e)}(\lambda_2 - \lambda_2 u) +1-\rho_1.
\end{eqnarray}

Let $T_{\alpha}^{(e)}$, $T_{\beta_i}^{(e)}$ and $T_{\beta}^{(e)}$ be r.v.s having the distributions $F_{\alpha}^{(e)}(x)$, $F_{\beta_i}^{(e)}(x)$ and $F_{\beta}^{(e)}(x)$, respectively.
From (\ref{Ka(u)}), we know
\begin{equation}\label{Ka-decomp}
K_a\stackrel{\rm d}{=}\left\{\begin{array}{ll}
0, &\mbox{ with probability }1-\rho_1,\\
N_{\lambda_2}(T_{\alpha}^{(e)}), &\mbox{ with probability }\rho_1.
\end{array}
\right.
\end{equation}

Next,  let $N_{\lambda,X_g}(T_{\beta}^{(e)})$ represent the number of the batched Poisson arrivals, with rate $\lambda$, and batch size $X_g$ within the time interval $(0,T_{\beta}^{(e)}]$. Then, by a similar conditioning argument as seen in (\ref{alpha-phi-0}), we have,
\begin{equation}\label{alpha-phi-0-b}
    E( z^{ N_{\lambda,X_g}(T_{\beta}^{(e)}) })=\int_0^{\infty}\sum_{n=0}^{\infty}(g(z))^n\frac {(\lambda x)^n} {n!} e^{-\lambda x}dF_{\beta}^{(e)}(x) = \beta^{(e)}(\lambda- \lambda g(z)),
\end{equation}
  where  $g(z)$ is the PGF of $X_g$. Now, it follows from (\ref{Kb-1}) that
\begin{equation}
    K_b(z) = \beta^{(e)}(\lambda- \lambda g(z)).
\end{equation}
Hence,
\begin{equation}\label{Kb-3}
K_b\stackrel{\rm d}{=}N_{\lambda,X_g}(T_{\beta}^{(e)}).
\end{equation}

Finally, from (\ref{Kc-1}), we have,
\begin{eqnarray}\label{Kc(u)}
K_c(u)&=& (1-\vartheta) \left[ 1 - \frac {1- g(u)} {1-u}\cdot \frac {1-\beta_2(\lambda-\lambda g(u))} {1- g(u)} \right]^{-1}\nonumber\\
&=& \frac {1-\vartheta} { 1 - \vartheta\cdot K_a(u)\cdot \beta_2^{(e)}(\lambda- \lambda g(u)) }\nonumber\\
&=&1-\vartheta +\vartheta\cdot\sum_{i=1}^{\infty}(1-\vartheta)\vartheta^{i-1}\big(K_a(u)\cdot \beta_2^{(e)}(\lambda- \lambda g(u))\big)^i,
\end{eqnarray}
where $\vartheta = \rho_2/(1-\rho_1)<1$.

A probabilistic interpretation for $K_c(u)$ is provided in the following remark for the convenience of future reference.
\begin{remark}\label{remark-Kc}
Let $\{X_c^{(i)}\}_{i=1}^{\infty}$ be a sequence of i.i.d. non-negative integer-valued r.v.s., each with the same PGF $K_a(z)\cdot \beta_2^{(e)}(\lambda- \lambda g(z))$, namely, $X_c^{(i)}\stackrel{\rm d}{=}K_a+N_{\lambda,X_g}(T_{\beta_2}^{(e)})$, where the two components are assumed to be independent. From (\ref{Kc(u)}), we know
\begin{equation}\label{Kc-decomp}
K_c\stackrel{\rm d}{=}\left\{\begin{array}{ll}
0, &\mbox{ with probability }1-\vartheta,\\
\sum_{i=1}^{J}X_c^{(i)}, &\mbox{ with probability }\vartheta,
\end{array}
\right.
\end{equation}
where $P(J=i)=(1-\vartheta)\vartheta^{i-1}$, $i\ge 1$, and $J$ is independent of $\{X_c^{(i)}\}_{i=1}^{\infty}$.
\end{remark}

Immediately from (\ref{K(u)-1}), we see that,
\begin{equation}\label{K-decomp}
K\stackrel{\rm d}{=}K_a+K_b+K_c,
\end{equation}
where $K_a$, $K_b$ and $K_c$ are assumed to be independent r.v.s.

\subsection{Stochastic decompositions on $(R_{11},R_{12})$ and $(R_{21},R_{22})$}

Recall $R_1(z_1,z_2)$ and $R_2(z_1,z_2)$ given in (\ref{P1-1}) and (\ref{P2-1}). Let
\begin{equation}
S_{\beta_i}(z_1,z_2)=\frac 1 {\beta_{i,1}}\cdot \frac {1-\beta_i(\lambda-\lambda_1 z_1 -\lambda_2 z_2)} {\lambda-\lambda_1 z_1 -\lambda_2 z_2}= \beta_i^{(e)}(\lambda -\lambda q z_1-\lambda p z_2),\quad i=1,2.\label{Ma-1}
\end{equation}
Simplifying (\ref{W-1}) gives us,
\begin{eqnarray}\label{P1-2}
W(z_1,z_2)
&=&\left (\lambda-\lambda g(z_2)\right ) \left[\beta_2(\lambda-\lambda_1 z_1-\lambda_2 z_2)-\beta_2(\lambda-\lambda g(z_2))\right]\nonumber \\
&&+ (\lambda_1 h(z_2) - \lambda_1 z_1  ) \left[z_2-\beta_2(\lambda-\lambda g(z_2))\right].
\end{eqnarray}
After substituting (\ref{P1-2}) into (\ref{P1-1}), we get
\begin{eqnarray}\label{R1-2}
R_1(z_1,z_2)&=&M_2(z_1,z_2)\cdot M_1(z_1,z_2)\cdot S_{\beta_1}(z_1,z_2)\cdot  R_0(z_2),
\end{eqnarray}
where
\begin{eqnarray}
M_1(z_1,z_2)&=&(1-\rho_1)\cdot \frac {h(z_2) - z_1} {\beta_1(\lambda-\lambda_1 z_1 - \lambda_2 z_2) - z_1},\label{Mb-1}\\
\\
M_2(z_1,z_2)&=&\frac {1-\rho} {\lambda_1(1-\rho_1)} \left[\frac {\left (\lambda-\lambda g(z_2)\right ) \left (\beta_2(\lambda-\lambda g(z_2))-\beta_2(\lambda-\lambda_1 z_1-\lambda_2 z_2)\right)} {(\beta_2(\lambda-\lambda g(z_2))-z_2) (h(z_2) - z_1)}+\lambda_1\right]. \nonumber \\
\label{Mc-1}
\end{eqnarray}
\\
Applying (\ref{Ka-1}) and (\ref{Kc-1}), we can rewrite (\ref{P2-1}) as
\begin{eqnarray}\label{R2-2}
R_2(z_1,z_2)= S_{\beta_2}(z_1,z_2)\cdot K_a(z_2)\cdot K_c(z_2)\cdot R_0(z_2).
\end{eqnarray}
Later, in Subsections \ref{subsection-Sbeta-i}, \ref{subsection3.4} and \ref{subsection-M2}, we will verify that $S_{\beta_i}(z_1,z_2)$ and $M_i(z_1,z_2)$ are the PGFs of two-dimensional r.v.s, denoted by $(S_{\beta_i,1},S_{\beta_i,2})$ and $(M_{i1},M_{i2})$, $i=1,2$, respectively.
Namely, $S_{\beta_i}(z_1,z_2)=E (z_1^{S_{\beta_i,1}}z_2^{S_{\beta_i,2}} )$ and $M_i(z_1,z_2)=E (z_1^{M_{i,1}}z_2^{M_{i,2}} )$, $i=1,2$.
Therefore, (\ref{R1-2}) and (\ref{R2-2}) imply that $(R_{11},R_{12})$ and $(R_{21},R_{22})$ can be decomposed into the sums of independent r.v.s.
Specifically,
\begin{eqnarray}
(R_{11},R_{12})&\stackrel{\rm d}{=}&(M_{21},M_{22})+(M_{11},M_{12})+(S_{\beta_1,1},S_{\beta_1,2})+(0, R_0), \label{R11-R12-decom}\\
(R_{21},R_{22})&\stackrel{\rm d}{=}&(S_{\beta_1,1},S_{\beta_1,2})+(0,K_a)+(0,K_c)+(0, R_0) \label{R21-R22-decom}.
\end{eqnarray}

\subsubsection{Probabilistic interpretation for the PGFs $S_{\beta_i}(z_1,z_2)$} \label{subsection-Sbeta-i}

For a probabilistic interpretation of the PGFs $S_{\beta_i}(z_1,z_2)$, $i=1,2$, let us introduce the following concept of splitting.
\begin{definition}
Let $N$ be a non-negative integer-valued r.v., and let $\{X_k\}_{k=1}^{\infty}$ be a sequence of i.i.d. Bernoulli r.v.s, which is independent of $N$, having the common $0$-$1$ distribution $P\{X_k=1\}=c$ and $P\{X_k=0\}=1-c$ with $0<c<1$. The two-dimensional r.v.
$(\sum_{k=1}^N X_k,N-\sum_{k=1}^N X_k)$, where $\sum_{1}^0\equiv 0$,
is called an independent $(c,1-c)$-splitting of $N$, denoted by $\splt(N;c,1-c)$.
\end{definition}

From the definition, it is easy to see that $(\sum_{k=1}^N X_k,N-\sum_{k=1}^N X_k)$ is independent of $(c,1-c)$-splitting of $N$, which is equivalent to
\begin{eqnarray}
    E \big (z_1^{\sum_{k=1}^N X_k} z_2^{N-\sum_{k=1}^N X_k} \big ) &=& \sum_{n=0}^{\infty} \left[\prod_{k=1}^n E(z_1^{X_k} z_2^{1- X_k})\right] P\{N=n\}\nonumber\\
     &=& \sum_{n=0}^{\infty} (cz_1+(1-c)z_2)^n P\{N=n\},\label{eqn:split}
\end{eqnarray}
where $\prod_{1}^0\equiv 1$.

A probabilistic interpretation of $S_{\beta_i}(z_1,z_2)$ is provided in terms of splitting in the following remark for the convenience of future reference.
\begin{remark}\label{M{a1},M{a2}}
For $i=1,2$, $(S_{\beta_i,1},S_{\beta_i,2})\stackrel{\rm d}{=} \splt(N_{\lambda}(T_{\beta_i}^{(e)});q,p)$,  since
\begin{equation}
S_{\beta_i}(z_1,z_2)=\beta_i^{(e)}(\lambda -\lambda_1 z_1-\lambda_2 z_2)=\int_0^{\infty}\sum_{n=0}^{\infty} (qz_1+pz_2)^n ((\lambda x)^n/n!) e^{-\lambda x}dF_{\beta_i}^{(e)}(x),
\end{equation}
which follows from (\ref{eqn:split}) by setting $N=N_{\lambda}(T_{\beta_i}^{(e)})$ and $c=q$.
\end{remark}

\subsubsection{Probabilistic interpretation for the PGF $M_1(z_1,z_2)$}\label{subsection3.4}

In this subsection, we prove that $M_1(z_1,z_2)$ is the PGF of a two-dimensional r.v. $(M_{11},M_{12})$. Let
\begin{eqnarray}
H_{\beta_1}(z_1,z_2)&=&\frac 1 {\rho_1}\cdot\frac {\beta_1(\lambda-\lambda_1 z_1-\lambda_2 z_2)-h(z_2)} {z_1-h(z_2)}\label{G2-decomp-H-0} \\
&=&\frac 1 {\rho_1}\cdot\frac {\beta_1(\lambda-\lambda_1 z_1-\lambda_2 z_2)-\beta_1(\lambda-\lambda_1 h(z_2)-\lambda_2 z_2)} {z_1-h(z_2)}\quad \mbox{(by (\ref{Falin-phi-eqn}))}.\label{G2-decomp-H}
\end{eqnarray}
It follows from (\ref{Mb-1}) and (\ref{G2-decomp-H-0}) that
\begin{equation}
M_1(z_1,z_2)=\frac {1-\rho_1} {1-\rho_1 H_{\beta_1}(z_1,z_2)}=\sum_{n=0}^{\infty} (1-\rho_1) \rho_1^{n} \left(H_{\beta_1}(z_1,z_2)\right)^n. \label{G2-decomp-eqn1}
\end{equation}
Clearly,
$(M_{11},M_{12})$ can be regarded as a random sum of two-dimensional r.v.s. provided that $H_{\beta_1}(z_1,z_2)$ is the PGF of a two-dimensional r.v.
To verify this, we will write (\ref{G2-decomp-H}) as a power series. Let $b_{\beta_1,k}=\int_0^{\infty}\frac{(\lambda t)^k} {k!}e^{-\lambda t}dF_{\beta_1}(t)$, $k\ge 1$.
Hence,
\begin{equation}
\beta_1(\lambda-\lambda q z_1-\lambda p z_2)=\int_0^{\infty}\sum_{k=0}^{\infty}\frac{(\lambda(q z_1+p z_2)t)^k} {k!}\cdot e^{-\lambda t}dF_{\beta_1}(t)=\beta_1(\lambda)+\sum_{k=1}^{\infty} b_{\beta_1,k} (q z_1+p z_2)^k. \label{H-lemma-1-proof-1}
\end{equation}
By (\ref{Falin-phi-eqn}) and (\ref{H-lemma-1-proof-1}),
\begin{equation}
h(z_2)=\beta_1(\lambda-\lambda q h(z_2)-\lambda p z_2)=\beta_1(\lambda)+\sum_{k=1}^{\infty} b_{\beta_1,k} (q h(z_2)+p z_2)^k. \label{H-lemma-1-proof-2}
\end{equation}
Substituting (\ref{H-lemma-1-proof-1}) and (\ref{H-lemma-1-proof-2}) into the numerator of the right-hand side of (\ref{G2-decomp-H}), we obtain
\begin{eqnarray}
H_{\beta_1}(z_1,z_2)
&=&\frac {1} {\rho_1}\cdot\sum_{k=1}^{\infty} b_{\beta_1,k} \left(\frac{(q z_1+p z_2)^k-(q h(z_2)+p z_2)^k} {z_1-h(z_2)}\right)\nonumber\\
&=&\frac {q} {\rho_1}\cdot\sum_{k=1}^{\infty} b_{\beta_1,k} \sum_{i=1}^{k}
(q z_1+p z_2)^{i-1}(q h(z_2)+p z_2)^{k-i}\nonumber\\
&=&\frac {q} {\rho_1}\cdot\sum_{k=1}^{\infty} k b_{\beta_1,k}\cdot D_{k}(z_1,z_2),\label{H-lemma-1-proof-3}
\end{eqnarray}
where
\begin{eqnarray}\label{Hj-lemma-1-eqn-1}
D_{k}(z_1,z_2)&=&\frac 1 k\sum_{i=1}^{k}
(q z_1+p z_2)^{i-1}(q h(z_2)+p z_2)^{k-i}.
\end{eqnarray}
Note that $q h(z_2)+p z_2=g(z_2)$ and $q z_1+p z_2$ are the PGFs of r.v.s (one or two-dimensional).
Hence, for $k\ge 1$, $D_{k}(z_1,z_2)$ is the PGF of a two-dimensional r.v., denoted by $(D_{k,1},D_{k,2})$.
In addition,
\begin{equation}\label{sum-bnbar}
\sum_{k=1}^{\infty} k b_{\beta_1,k}=\sum_{k=1}^{\infty} \int_0^{\infty}\frac{(\lambda t)^k} {(k-1)!}e^{-\lambda t}dF_{\beta_1}(t)=\lambda\int_0^{\infty}tdF_{\beta_1}(t)=\rho_1/q.
\end{equation}
Namely, $(q/\rho_1)\sum_{k=1}^{\infty}kb_{\beta_1,k}=1$, which together with (\ref{H-lemma-1-proof-3}) implies that $H_{\beta_1}(z_1,z_2)$ is the PGF of a   two-dimensional r.v., denoted by $(H_{\beta_1,1},H_{\beta_1,2})$. Namely, $H_{\beta_1}(z_1,z_2)=E (z_1^{H_{\beta_1,1}}z_2^{H_{\beta_1,2}})$.
The above argument is summarized in the following remarks.
\begin{remark}
Suppose that $\{(Y_{n,1},Y_{n,2})\}_{n=1}^{\infty}$ is a sequence of independent two-dimensional r.v.s, each with a common PGF $q z_1+p z_2$,  $\{Z_m\}_{m=1}^{\infty}$ is a sequence of independent r.v.s, each with a common PGF $g(z_2)$, and the two sequences are independent.
It follows from (\ref{Hj-lemma-1-eqn-1}) that for $k\ge 1$,
\begin{equation}
(D_{k,1},D_{k,2})\overset{\footnotesize\mbox{d}}{=}
\sum_{n=1}^{i-1}(Y_{n,1},Y_{n,2})+\sum_{m=1}^{k-i}(0,Z_{m}) \mbox{ with probability }1/k,\quad i=1,2,\cdots,k.
\end{equation}
\end{remark}
\begin{remark}\label{remark-H{11}H{12}}
It follows from (\ref{H-lemma-1-proof-3}) that
\begin{equation}
(H_{\beta_1,1},H_{\beta_1,2})\overset{\footnotesize\mbox{d}}{=}
(D_{k,1},D_{k,2}) \mbox{ with probability }(q/\rho_1)kb_{\beta_1,k},\quad k\ge 1.
\end{equation}
\end{remark}
\begin{remark}\label{remark-M{b1}M{b2}}
It follows from (\ref{G2-decomp-eqn1}) that $(M_{11},M_{12})$ is a random sum of i.i.d. two-dimensional r.v.s $(H_{\beta_1,1}^{(i)},H_{\beta_1,2}^{(i)})$, $i\ge 1$, each with the same PGF $H_{\beta_1}(z_1,z_2)$, and precisely,
\begin{equation}\label{remark-M{b1}M{b2}-eqn1}
(M_{11},M_{12})\stackrel{\rm d}{=}\left\{\begin{array}{ll}
0, &\mbox{ with probability }1-\rho_1,\\
\sum_{i=1}^{J}(H_{\beta_1,1}^{(i)},H_{\beta_1,2}^{(i)}) &\mbox{ with probability }\rho_1,
\end{array}
\right.
\end{equation}
where $P(J=i)=(1-\rho_1)\rho_1^{i-1}$, $i\ge 1$ and $J$ is independent of $(H_{\beta_1,1}^{(i)},H_{\beta_1,2}^{(i)})$, $i\ge 1$.
\end{remark}

\subsubsection{Probabilistic interpretation for the PGF $M_2(z_1,z_2)$}\label{subsection-M2}

In this subsection, we prove that $M_2(z_1,z_2)$ is the PGF of a two-dimensional r.v. $(M_{21},M_{22})$.
Let
\begin{eqnarray}
H_{\beta_2}(z_1,z_2)
&=&\frac p {q\rho_2}\cdot\frac {\beta_2(\lambda-\lambda_1 z_1-\lambda_2 z_2)-\beta_2(\lambda-\lambda_1 h(z_2)-\lambda_2 z_2)} {z_1-h(z_2)}.\label{H{beta_2}-1}
\end{eqnarray}
Using (\ref{H{beta_2}-1}), (\ref{Ka-1}) and (\ref{Kc-1}), we can rewrite (\ref{Mc-1}) as follows:
\begin{equation}\label{M_2(z_1,z_2)-34}
M_2(z_1,z_2)=\vartheta\cdot H_{\beta_2}(z_1,z_2)\cdot K_a(z_2)\cdot K_c(z_2)+ 1-\vartheta.
\end{equation}
It can be shown that $H_{\beta_2}(z_1,z_2)$ is the PGF of a two-dimensional r.v., denoted by $(H_{\beta_2,1},H_{\beta_2,2})$.
Namely, $H_{\beta_2}(z_1,z_2)=E (z_1^{H_{\beta_2,1}}z_2^{H_{\beta_2,2}})$.
The proof is similar to that for $H_{\beta_1}(z_1,z_2)$ in Subsection~\ref{subsection3.4}, details of which are omitted here.
\\

Let $b_{\beta_2,k}=\int_0^{\infty}\frac{(\lambda t)^k} {k!}e^{-\lambda t}dF_{\beta_2}(t)$, $k\ge 1$.
It follows from (\ref{H{beta_2}-1}) that
\begin{equation}\label{H{beta_2}-3}
H_{\beta_2}(z_1,z_2)
=\frac {p} {\rho_2}\cdot\sum_{k=1}^{\infty} k b_{\beta_2,k}\cdot D_{k}(z_1,z_2).
\end{equation}
Similar to (\ref{sum-bnbar}), we can verify that $(p /\rho_2)\sum_{k=1}^{\infty}k b_{\beta_2,k}=1$, which together with (\ref{H{beta_2}-3}) implies that $H_{\beta_2}(z_1,z_2)$ is the PGF of a two-dimensional r.v.
The above argument leads to the following two remarks.
\begin{remark}\label{remark-H{21}H{22}}
It follows from (\ref{H{beta_2}-3}) that
\begin{equation}
(H_{\beta_2,1},H_{\beta_2,2})\overset{\footnotesize\mbox{d}}{=} (D_{k,1},D_{k,2}) \mbox{ with probability }(p/\rho_2)kb_{\beta_2,k}.
\end{equation}
\end{remark}
\begin{remark}\label{remark-M{21}M{22}}
It follows from (\ref{M_2(z_1,z_2)-34}) that
\begin{equation}\label{def-T-theta}
(M_{21},M_{22})\stackrel{\rm d}{=}\left\{\begin{array}{ll}
0, &\mbox{ with probability }1-\vartheta,\\
(H_{\beta_2,1},H_{\beta_2,2})+(0, K_a)+(0, K_c), &\mbox{ with probability }\vartheta.
\end{array}
\right.
\end{equation}
\end{remark}

\section{Tail Asymptotics} \label{sec:4}

In this section, we study the asymptotic behaviour of the tail probabilities
$P\{R_0>j\}$ and $P\{R_{ik}>j\}$, $i,k=1,2$, as $j\to\infty$.

Applying Karamata's theorem (e.g., p.28 in \cite{Bingham:1989}), and using Assumption A1 and Lemma~\ref{Lemma 2.1}, respectively, gives,
as $t\to\infty$,
\begin{eqnarray}
P\{T_{\beta_1}^{(e)}>t\}&\sim& \frac {\lambda_1} {\rho_1(a_1-1)} \cdot t^{-a_1+1} L(t),\label{T-beta-i}\\
P\{T_{\alpha}^{(e)}>t\}&\sim& \frac {1} {\alpha_1 (a_1-1)(1-\rho_1)^{a_1+1}} \cdot t^{-a_1+1} L(t).\label{T-alpha}
\end{eqnarray}

Applying Proposition 8.5 (p.181 in \cite{Grandell:1997}) to the density $\bar{F}_{\beta_2}(t)/\beta_{2,1}$ and using Assumption A2, gives,
as $t\to\infty$,
\begin{eqnarray}\label{T-beta-2}
P\{T_{\beta_2}^{(e)}>t\}&\sim&\left\{\begin{array}{ll}
\frac {\lambda_2} {\rho_2 (a_2-1) } \cdot t^{-a_2+1} L(t), &\mbox{ if }r=0,\\
\frac {\lambda_2} {\rho_2 r}\cdot e^{-rt} t^{-a_2} L(t), &\mbox{ if }r>0.
\end{array}
\right.
\end{eqnarray}
Furthermore, since $F_{\beta}(x)=q F_{\beta_1}(x)+p F_{\beta_2}(x)$ and based on Assumptions~ A1 and A2, we have, $P\{T_{\beta}>t\}=qP\{T_{\beta_1}>t\} +p P\{T_{\beta_2}>t\}\sim q t^{-a_1} L(t)$ as $t\to\infty$, from which Karamata's theorem implies that
\begin{eqnarray}
P\{T_{\beta}^{(e)}>t\}&\sim& \frac {\lambda_1} {\rho(a_1-1) } \cdot t^{-a_1+1} L(t).\label{T-beta}
\end{eqnarray}

\subsection{Asymptotic tail probability of the r.v. $R_0$}\label{R0+Mc}

Recall (\ref{D0-2}), which closely relates the PGF of $R_0$ to the PGF of $K$. For this reason, we first study the tail probability for $K$, which can be regarded as a sum of independent r.v.s $K_a$, $K_b$ and $K_c$ (refer to (\ref{K-decomp})).
By (\ref{Ka-decomp}), (\ref{T-alpha}) and applying Lemma~\ref{lemma-Asmussen-poisson}, we have,
\begin{eqnarray}\label{Ka>j}
P\{K_a>j\}&=&\rho_1 P\{N_{\lambda_2}(T_{\alpha}^{(e)})>j\}
\sim \frac {\lambda_1\lambda_2^{a_1-1}}{(a_1-1) (1-\rho_1)^{a_1}} \cdot j^{-a_1+1} L(j),\quad j\to\infty.
\end{eqnarray}

Recall (\ref{Kb-3}), $K_b=N_{\lambda, X_g}(T_{\beta}^{(e)})\stackrel{\rm d}{=}\sum_{i=1}^{N_{\lambda}(T_{\beta}^{(e)})} X_g^{(i)}$, where $X_g^{(i)}$ has the common distribution $X_g$. By (\ref{G-def}), and then applying Lemma~\ref{lemma-Asmussen-poisson} and using Lemma~\ref{Lemma 2.1}, we know that
\[
    P\{X_g>j\} \sim  qP\{N_{\lambda_2}(T_{\alpha})>j\}\sim q\lambda_2^{a_1}(1-\rho_1)^{-a_1-1} \cdot j^{-a_1} L(j).
\]
Similarly, applying Lemma~\ref{lemma-Asmussen-poisson} and using (\ref{T-beta}), we have,
\[
    P\{N_{\lambda}(T_{\beta}^{(e)})>j\} \sim \frac {q\lambda^{a_1-1}} {(q\beta_{1,1}+p\beta_{2,1})  (a_1-1) } \cdot  j^{-a_1+1} L(j).
\]
Based on which, by (\ref{G-exp}) and applying Lemma~\ref{Lemma 4.5}, we have,
\begin{eqnarray}
P\{K_b>j\}&\sim&\frac {\lambda_1 \lambda_2^{a_1-1}} {\rho (a_1-1)(1-\rho_1)^{a_1-1}} \cdot j^{-a_1+1} L(j),
\quad j\to\infty.\label{Kb>j}
\end{eqnarray}

Next, we study $P\{K_c>j\}$. By Remark~\ref{remark-Kc}, we know that
$P\{K_c>j\}=\vartheta P\{\sum_{i=1}^{J}X_c^{(i)}>j\}$,
where $P(J=i)=(1-\vartheta)\vartheta^{i-1}$, $i\ge 1$, and $X_c^{(i)}$ has the same distribution as $X_c=K_a+N_{\lambda,X_g}(T_{\beta_2}^{(e)})$.
Note that $N_{\lambda,X_g}(T_{\beta_2}^{(e)})\stackrel{\rm d}{=}\sum_{i=1}^{N_{\lambda}(T_{\beta_2}^{(e)})} X_g^{(i)}$, where $X_g^{(i)}$ has the common tail probability $P\{X_g>j\}\sim Const\cdot j^{-a_1} L(j)$ and $P\{N_{\lambda}(T_{\beta_2}^{(e)})>j\}\sim Const\cdot j^{-a_2+1}L(j)$, where the symbol ``$Const$" stands for a constant, and such a symbol will be used throughout the paper. Therefore, by applying Lemma~\ref{Lemma 4.5} (and noticing that $a_2 > a_1$ if $r=0$ in Assumption~A2),
\begin{eqnarray}
P\{N_{\lambda, X_g}(T_{\beta_2}^{(e)})>j\} \sim Const\cdot\max (j^{-a_2+1}L(j),j^{-a_1}L(j))=o(1)\cdot j^{-a_1+1} L(j).\label{N{lambda,Xg}}
\end{eqnarray}
By (\ref{Ka>j}), (\ref{N{lambda,Xg}}), applying Lemma~\ref{Embrechts-compound-geo} and Lemma~\ref{Lemma 3.4new}, we have, as $j\to\infty$,
\begin{eqnarray}
P\{K_c>j\}\sim\frac {\vartheta} {1-\vartheta}P\{X_c>j\}= \frac {\vartheta} {1-\vartheta} P\{K_a+N_{\lambda, X_g}(T_{\beta_2}^{(e)})>j\}\sim\frac {\vartheta} {1-\vartheta} P\{K_a>j\},
\end{eqnarray}
which, together with (\ref{Ka>j}), (\ref{Kb>j}) and (\ref{K-decomp}), leads to
\begin{eqnarray}
P\{K_a+K_c>j\}&\sim& \frac {\lambda_1\lambda_2^{a_1-1}}{(1-\rho) (a_1-1) (1-\rho_1)^{a_1-1}} \cdot j^{-a_1+1} L(j),\quad j\to\infty,\label{Ka+Kc-asym}\\
P\{K>j\}&\sim& \frac {\lambda_1\lambda_2^{a_1-1}}{\rho(1-\rho)(a_1-1) (1-\rho_1)^{a_1-1}} \cdot j^{-a_1+1} L(j),\quad j\to\infty.
\end{eqnarray}

By (\ref{D0-2}), the PGF $R_0(z)$ is expressed in terms of the PGF $K(z)$.
Therefore, the tail probability of $R_0$ is determined by the tail probability of $K$. The following asymptotic result is a straightforward application
of Theorem 5.1 in~\cite{Liu-Min-Zhao:2017}.
\begin{eqnarray}
P\{R_0>j\}&\sim&\frac {a_1-1} {a_1} \psi \cdot \frac {\lambda_1\lambda_2^{a_1-1}}{\rho(1-\rho)(a_1-1) (1-\rho_1)^{a_1-1}} \cdot j^{-a_1} L(j),\label{P{R_0>j}sim-0}
\end{eqnarray}
where $\psi$ is given in (\ref{psi}).
Recall the definition of $R_{0}$ in Section 2. The above discussion is summarized in the following theorem.
\begin{theorem} \label{the:PO-0}
As $j\to\infty$,
\begin{eqnarray}
P\{R_{orb}>j|I_{ser}=0\}=P\{R_{0}>j\}&\sim&\frac {\lambda_1 \lambda_2^{a_1}} {a_1\mu(1-\rho)^2 (1-\rho_1)^{a_1-1}}\cdot j^{-a_1}L(j).\label{P{R_0>j}sim}
\end{eqnarray}
\end{theorem}

\subsection{Asymptotic tail probabilities of the r.v.s $R_{11}$, $R_{12}$, $R_{21}$ and $R_{22}$}

Recalling (\ref{R11-R12-decom}) and (\ref{R21-R22-decom}), we immediately have
\begin{eqnarray}
R_{11} &\stackrel{\rm d}{=}& M_{21}+M_{11}+S_{\beta_1,1},\label{R11-decom}\\
R_{12} &\stackrel{\rm d}{=}& M_{22}+M_{12}+S_{\beta_1,2}+R_0,\label{R12-decom}\\
R_{21} &\stackrel{\rm d}{=}& S_{\beta_2,1},\label{R21-decom}\\
R_{22} &\stackrel{\rm d}{=}& S_{\beta_2,2}+ K_a+ K_c +R_0,\label{R22-decom}
\end{eqnarray}
where all of the r.v.s on the right hand side in each of (\ref{R11-decom})--(\ref{R22-decom}) are independent.

In the previous sections, the asymptotic behaviour of the tail probabilities for the r.v.s $R_0$ and $ K_a+ K_c$ have already been obtained in (\ref{P{R_0>j}sim}) and (\ref{Ka+Kc-asym}), respectively.
In the following, we will focus on the tail probabilities of the r.v.s $S_{\beta_i,k}$ and $M_{ik}$ for $i,k=1,2$.

Recall Remark~\ref{M{a1},M{a2}}, $S_{\beta_1,1}\stackrel{\rm d}{=}N_{\lambda_1}(T_{\beta_1}^{(e)})$ and $S_{\beta_1,2}\stackrel{\rm d}{=}N_{\lambda_2}(T_{\beta_1}^{(e)})$, $S_{\beta_2,1}\stackrel{\rm d}{=}N_{\lambda_1}(T_{\beta_2}^{(e)})$ and $S_{\beta_2,2}\stackrel{\rm d}{=}N_{\lambda_2}(T_{\beta_2}^{(e)})$. By (\ref{T-beta-i}) and applying Lemma~\ref{lemma-Asmussen-poisson}, we obtain
\begin{eqnarray}
P\{S_{\beta_1,1}>j\}&\sim& \frac {\lambda_1^{a_1}} {\rho_1 (a_1-1) } \cdot j^{-a_1+1} L(j)\label{Sbeta11-tail}, \\
P\{S_{\beta_1,2}>j\}&\sim& \frac {\lambda_1\lambda_2^{a_1-1}} {\rho_1 (a_1-1) } \cdot j^{-a_1+1} L(j)\label{Sbeta12-tail}.
\end{eqnarray}
By (\ref{T-beta-2}) and applying Lemma~\ref{lemma-Asmussen-poisson} and Lemma~\ref{lemma-Grandell-poisson}, we obtain
\begin{eqnarray}
P\{S_{\beta_2,1}>j\}&\sim&\left\{\begin{array}{ll}
\frac {\lambda_2\lambda_1^{a_2-1}} {\rho_2 (a_2-1)}\cdot j^{-a_2+1} L(j), &\mbox{ if }r=0,\\
\frac {\lambda_2\lambda_1 (\lambda_1+r)^{a_2-1}} {\rho_2 r} \cdot \left(\frac {\lambda_1} {\lambda_1+r}\right)^j j^{-a_2} L(j), &\mbox{ if }r>0.
\end{array}
\right. \label{Sbeta21-tail} \\
P\{S_{\beta_2,2}>j\}&\sim&\left\{\begin{array}{ll}
\frac {\lambda_2^{a_2}} {\rho_2 (a_2-1) }\cdot j^{-a_2+1} L(j), &\mbox{ if }r=0,\\
\frac{\lambda_2^2(\lambda_2+r)^{a_2-1}} {\rho_2 r}\cdot  \left(\frac {\lambda_2} {\lambda_2+r}\right)^j j^{-a_2} L(j), &\mbox{ if }r>0.
\end{array}
\right.\label{Sbeta22-tail}
\end{eqnarray}
\\
Next, we will study the asymptotic tail probabilities of the r.v.s $M_{ik},\ i,k=1,2$.
By Remark \ref{remark-M{b1}M{b2}} and Remark \ref{remark-M{21}M{22}}, we know that
\begin{eqnarray}
P\{M_{1k}>j\}&=&\rho_1 P\Big\{\sum_{i=1}^{J}H_{\beta_1,k}^{(i)}>j\Big\},\quad k=1,2,\label{M1k>j}\\
P\{M_{21}>j\}&=&\vartheta P \{H_{\beta_2,1} >j \},\label{M21>j}\\
P\{M_{22}>j\}&=&\vartheta P\{H_{\beta_2,2}+K_a+K_c >j \}.\label{M22>j}
\end{eqnarray}
To proceed further, we need to study the tail probabilities of the r.v.s $H_{\beta_i,k}$ for $i,k=1,2$.

\subsubsection{Asymptotic tail probabilities of the r.v.s $H_{\beta_1,1}$ and $H_{\beta_2,1}$}

Taking $z_2\to 1$ in (\ref{G2-decomp-H}) and (\ref{H{beta_2}-1}), we can write
\begin{eqnarray}
E(z_1^{H_{\beta_1,1}})&=&H_{\beta_1}(z_1,1)=\frac 1 {\rho_1}\cdot\frac {\beta_1(\lambda_1-\lambda_1 z_1)-1} {z_1-1}=\beta_1^{(e)}(\lambda_1-\lambda_1 z_1),\\
E(z_1^{H_{\beta_2,1}})&=&H_{\beta_2}(z_1,1)=\frac p {q\rho_2}\cdot\frac {\beta_2(\lambda_1-\lambda_1 z_1)-1} {z_1-1}=\beta_2^{(e)}(\lambda_1-\lambda_1 z_1).
\end{eqnarray}
Therefore, $H_{\beta_i,1}\stackrel{\rm d}{=}N_{\lambda_1}(T_{\beta_i}^{(e)})\stackrel{\rm d}{=}S_{\beta_i,1}$, $i=1,2$, and
\begin{equation}
P\{H_{\beta_i,1}>j\}=P\{S_{\beta_i,1}>j\},\quad i=1,2,\label{Sbeta-i1-tail}
\end{equation}
whose asymptotic tails are present in (\ref{Sbeta11-tail}) and (\ref{Sbeta21-tail}), respectively.

\subsubsection{Asymptotic tail probability of the r.v. $H_{\beta_1,2}$}\label{subsectionH12}

Unlike the other r.v.s discussed early, more efforts are required for the asymptotic tail behaviour for $H_{\beta_1,2}$, which will be presented in Proposition~\ref{lemma-Hbeta12>j}.
Before doing that, we first present a nice bound on the tail probability of $H_{\beta_1,2}$,
which is very illustrative for an intuitive understanding of the tail property for $H_{\beta_1,2}$.

Taking $z_1\to 1$ in (\ref{Hj-lemma-1-eqn-1}) and (\ref{H-lemma-1-proof-3}), we have,
\begin{eqnarray}
E(z_2^{D_{k,2}})=D_{k}(1,z_2)&=&\frac 1 k\sum_{i=1}^{k}(q +p z_2)^{i-1}(q h(z_2)+p z_2)^{k-i}, \label{Dk(1,z)}\\
E(z_2^{H_{\beta_1,2}})=H_{\beta_1}(1,z_2)&=&\frac {q} {\rho_1}\cdot\sum_{k=1}^{\infty} k b_{\beta_1,k}\cdot D_{k}(1,z_2).\label{Hbeta1(1,z)}
\end{eqnarray}
It follows from (\ref{Dk(1,z)}) that for $k\ge 1$,
$$
D_{k,2}\stackrel{\rm d}{=}\sum_{n=1}^{i-1}Y_n+\sum_{n=1}^{k-i}Z_n\quad\mbox{with probability } 1/k\quad\mbox{for }i=1,2,\cdots,k,
$$
where $\{Y_n\}_{n=1}^{\infty}$ and $\{Z_n\}_{n=1}^{\infty}$ are sequences of independent r.v.s that are independent of each other, with
$Y_n$ and $Z_n$ having PGFs $q +p z_2$ and $q h(z_2)+p z_2$, respectively.

We say that $Y$ is stochastically smaller than $Z$, written as $Y\le_{st}Z$, if $P\{Y>t\}\le P\{Z>t\}$ for all $t$. It is easy to see that $Y_{n_1}\le_{st}Z_{n_2}$ for all $n_1,n_2\ge 1$. Define
$$
D^L_{k,2}\stackrel{\rm d}{=}\sum_{n=1}^{k-1}Y_n\quad\mbox{ and }\quad D^U_{k,2}\stackrel{\rm d}{=}\sum_{n=1}^{k-1}Z_n.
$$
Then, by Theorem 1.2.17 (p.7 in \cite{Muller-Stoyan:2002}),
\begin{equation}
D^L_{k,2}\le_{st}D_{k,2}\le_{st}D^U_{k,2}.\label{DDD}
\end{equation}
Furthermore, it follows from (\ref{Hbeta1(1,z)}) that
$H_{\beta_1,2}\stackrel{\rm d}{=}D_{k,2}$, with probability $(q/\rho_1)k b_{\beta_1,k}$, for $k\ge 1$.
\\
\\
Now define the r.v.s $H^L_{\beta_1,2}$ and $H^U_{\beta_1,2}$ as follows:
\begin{eqnarray}
H^L_{\beta_1,2}&\stackrel{\rm d}{=}&D^L_{k,2}\quad\mbox{with probability } (q/\rho_1)k b_{\beta_1,k}\ \mbox{for }k\ge 1,\nonumber\\
H^U_{\beta_1,2}&\stackrel{\rm d}{=}&D^U_{k,2}\quad\mbox{with probability } (q/\rho_1)k b_{\beta_1,k}\ \mbox{for }k\ge 1.\nonumber
\end{eqnarray}
Then, by (\ref{DDD}),
\begin{equation}\label{HU-H-HL}
H_{\beta_1,2}^{L}\le_{st}H_{\beta_1,2}\le_{st}H_{\beta_1,2}^{U}.
\end{equation}
Note that $H_{\beta_1,2}^{L}$ and $H_{\beta_1,2}^{U}$ have the following PGFs:
\begin{eqnarray}
E(z_2^{H_{\beta_1,2}^{L}})&=&\frac {q} {\rho_1}\cdot\sum_{k=1}^{\infty} k b_{\beta_1,k}\cdot E(z_2^{D_{k,2}^{L}})
=\frac {q} {\rho_1}\cdot\sum_{k=1}^{\infty} k b_{\beta_1,k}\cdot (q +p z_2)^{k-1},\label{HbetaL}\\
E(z_2^{H_{\beta_1,2}^{U}})&=&\frac {q} {\rho_1}\cdot\sum_{k=1}^{\infty} k b_{\beta_1,k}\cdot E(z_2^{D_{k,2}^{U}})
=\frac {q} {\rho_1}\cdot\sum_{k=1}^{\infty} k b_{\beta_1,k}\cdot (qh(z_2) +p z_2)^{k-1}.\label{HbetaU}
\end{eqnarray}

Next, we will study the asymptotic behaviour of $P\{H_{\beta_1,2}^{L}>j\}$ and $P\{H_{\beta_1,2}^{U}>j\}$, respectively.
Let $N$ be a r.v. with probability distribution $P\{N=k\}=(q/\rho_1) k b_{\beta_1,k}$, $k\ge 1$. Therefore, (\ref{HbetaL}) and (\ref{HbetaU}) can be written as
\[
H_{\beta_1,2}^{L}\stackrel{\rm d}{=}\sum_{k=1}^{N-1} Y_k
\quad\mbox{and}\quad
H_{\beta_1,2}^{U}\stackrel{\rm d}{=}\sum_{k=1}^{N-1} Z_k,
\]
where $N$ is independent of both $Z_k$ and $Y_k$, $k\ge 1$.
\\
\\
Then, it is immediately clear that,
\begin{equation}
P\{N>m\}=(q/\rho_1)\sum_{k=m+1}^{\infty} k b_{\beta_1,k}=(q/\rho_1)\left[m\overline{b}_{\beta_1,m+1}+\sum_{k=m+1}^{\infty} \overline{b}_{\beta_1,k}\right],\label{P{N>m}}
\end{equation}
where $\overline{b}_{\beta_1,k}=\sum_{n=k}^{\infty} b_{\beta_1,n}$.
\\
\\
Using the definition of $b_{\beta_1,n}$ in Section~\ref{subsection3.4}, and by applying Lemma~\ref{lemma-Asmussen-poisson}, we know $\overline{b}_{\beta_1,k}=P\{N_{\lambda}(T_{\beta_1})>k-1\}\sim \lambda^{a_1}k^{-a_1}L(k)$ as $k\to\infty$, which, together with Proposition~1.5.10 in \cite{Bingham:1989}, implies that
\begin{eqnarray}
P\{N>m\}&\sim& \frac {a_1 q\lambda^{a_1}}{\rho_1(a_1-1)}m^{-a_1+1}L(m)\quad{as}\ m\to\infty.
\end{eqnarray}
Recall the following three facts: (i) $Y_k$ is a $0-1$ r.v., which implies that $P\{Y_k>j\}\to 0$ as $j\to\infty$; (ii) $Z_k$ has the same probability distribution as $X_g$ defined in (\ref{G-def}), which implies that $P\{Z_k>j\}=P \{X_g>j\}\sim Const\cdot j^{-a_1} L(j)$ as $j\to\infty$; and (iii)
$E(Y_k)=p$ and $E(Z_k)=E(X_g)=p/(1-\rho_1)<\infty$ given in (\ref{G-exp}). Then, by Lemma~\ref{Lemma 4.5}, we know
\begin{eqnarray}
P\{H_{\beta_1,2}^{L}>j\}&\sim& a_1\cdot\frac{\lambda_1 \lambda_2^{a_1-1}} {\rho_1(a_1-1)}\cdot j^{-a_1+1}L(j)\label{HbetaL-34}\\
P\{H_{\beta_1,2}^{U}>j\}&\sim&\frac {a_1}{(1-\rho_1)^{a_1-1}}\cdot\frac{\lambda_1 \lambda_2^{a_1-1}} {\rho_1(a_1-1)}\cdot j^{-a_1+1}L(j).\label{HbetaU-35}
\end{eqnarray}

\begin{remark}\label{BoundAsymp}
It follows from (\ref{HU-H-HL}) that $P\{H_{\beta_1,2}^{L}>j\}\le P\{H_{\beta_1,2}>j\}\le P\{H_{\beta_1,2}^{U}>j\}$,
whereas the asymptotic properties of $P\{H_{\beta_1,2}^{L}>j\}$ and $P\{H_{\beta_1,2}^{U}>j\}$ are given in (\ref{HbetaL-34}) and (\ref{HbetaU-35}), respectively. This suggests that $P\{H_{\beta_1,2}>j\}\sim c\cdot\frac{\lambda_1 \lambda_2^{a_1-1}} {\rho_1(a_1-1)}\cdot j^{-a_1+1}L(j)$ as $j\to\infty$ for some constant $c\in \left(a_1, a_1/(1-\rho_1)^{a_1-1}\right)$. In the following proposition (Proposition~\ref{lemma-Hbeta12>j}), we will verify that this assertion is true.
\end{remark}

\begin{proposition}\label{lemma-Hbeta12>j}
As $j\to\infty$,
\begin{equation}
P\{H_{\beta_1,2}>j\}\sim\frac{1-\rho_1} {\rho_1} \left[\frac{1} {(1-\rho_1)^{a_1}}-1\right]\cdot\frac{\lambda_1 \lambda_2^{a_1-1}} {\rho_1(a_1-1)}\cdot j^{-a_1+1}L(j).\label{lemma-Hbeta12>j-formula}
\end{equation}
\end{proposition}

To prove this proposition, we need the following two lemmas (Lemma~\ref{lem:gamma is LST} and Lemma~\ref{proposition-A1}).
Setting $z_1=1$ in (\ref{G2-decomp-H-0}) and noting $h(z_2)=\alpha(\lambda_2-\lambda_2 z_2)$, we obtain
\begin{eqnarray}\label{PGF-Hbeta12}
E(z_2^{H_{\beta_1,2}})=H_{\beta_1}(1,z_2)
&=&\frac 1 {\rho_1}\cdot\frac {\beta_1(\lambda_2-\lambda_2 z_2)-\alpha(\lambda_2-\lambda_2 z_2)} {1-\alpha(\lambda_2-\lambda_2 z_2)}
=\gamma(\lambda_2-\lambda_2 z_2),
\end{eqnarray}
where
\begin{eqnarray}
\gamma(s)&=&\frac 1 {\rho_1}\cdot\frac {\beta_1(s)-\alpha(s)} {1-\alpha(s)}.\label{Apdx-1}
\end{eqnarray}

\begin{lemma} \label{lem:gamma is LST}
$\gamma(s)$ is the LST of a probability distribution on $[0,\infty)$.
\end{lemma}

\proof
By Theorem~1 in Feller (1991)~\cite{Feller1971} (see p.439),
it is true iff $\gamma(0)=1$ and $\gamma(s)$ is completely monotone, i.e., $\gamma(s)$ possesses derivatives of all orders such that $(-1)^{n}\frac {d^n}{ds^n}\gamma(s)\ge 0$ for $s> 0$, $n\ge 0$. It is easy to check by (\ref{Apdx-1}) that $\tau(0)=1$. Next, we verify that $\gamma(s)$ is completely monotone by using Criterion A.1 and Criterion A.2 in the appendix.

\noindent{\it Fact 1.} Take $\vartheta_3(s)=1/s$ and $\vartheta_4(s)=1-\alpha(s)\ge 0$ for $s> 0$. Because
$(-1)^n \frac {d^n\vartheta_3(s)} {ds^n}=\frac {n!} {s^{n+1}}> 0$ for $s> 0$ and $(-1)^n \frac {d^n\vartheta'_4(s)} {ds^n}=(-1)^{n+1} \frac {d^{n+1}\alpha(s)} {ds^{n+1}}\ge 0$ for $s> 0$, both $\vartheta_3(s)$ and $\vartheta'_4(s)$ are completely monotone. By Criterion A.2, we know that $1/(1-\alpha(s))$ is completely monotone.

\noindent{\it Fact 2.} It can be shown that $\beta_1(s)-\alpha(s)$ is completely monotone, i.e., $(-1)^n(\beta_1^{(n)}(s)-\alpha^{(n)}(s))\ge 0$ for $s> 0$, $n\ge 0$, where $\beta_1^{(n)}(s)$ and $\alpha^{(n)}(s)$ represent the $n$th derivative of $\beta_1^{(n)}(s)$ and $\alpha^{(n)}(s)$, respectively.
Let us proceed with using mathematical induction on $n\ge 0$.
Clearly, it is true for $n=0$ because $\beta_1(s)\ge \beta_1(s+\lambda_1- \lambda_1 \alpha(s))=\alpha(s)$ (by (\ref{busy-eqn-alpha})).
Now, let us make the induction hypothesis that $(-1)^{k}\kappa^{(k)}(s)\ge 0$ for $s> 0$ and all $k=0,1,\ldots, n$.
Then, by the mean value theorem, for $n\ge 0$, there exists some $c_n\in (0, 1)$ such that
\begin{eqnarray}
\beta_1^{(n)}(s)-\alpha^{(n)}(s)&=&\beta_1^{(n)}(s)-\beta_1^{(n)}(s+\lambda_1- \lambda_1 \alpha(s))\nonumber\\
&=&-\beta_1^{(n+1)}(s+c_n\cdot(\lambda_1- \lambda_1 \alpha(s)))\cdot(\lambda_1- \lambda_1 \alpha(s)).\label{Apdx-beta-alpha}
\end{eqnarray}
Note that $(-1)^n\beta_1^{(n)}(s)\ge 0$ for $s> 0$, $n\ge 0$.  The result (\ref{Apdx-beta-alpha}), together with the induction hypothesis, completes the proof for $k=n+1$.

By (\ref{Apdx-1}), Facts 1 and 2, and applying Criterion A.1, we know that $\gamma(s)$ is completely monotone. Therefore, it is the LST of a probability distribution.
\pend

\begin{remark}\label{gamma-compound-possion}
Let $T_{\gamma}$ be a r.v. whose the probability distribution has the LST $\gamma(s)$.
Then the expression $Ez_2^{H_{\beta_{1,2}}}=\gamma(\lambda_2-\lambda_2 z_2)$,  in (\ref{PGF-Hbeta12}), implies that $H_{\beta_{1,2}}\stackrel{\rm d}{=}N_{\lambda_2}(T_{\gamma})$.
\end{remark}

\begin{lemma}\label{proposition-A1}
As $t\to\infty$,
\begin{eqnarray}
P\{T_{\gamma}>t\}&\sim&\frac{1-\rho_1} {\rho_1} \left[\frac{1} {(1-\rho_1)^{a_1}}-1\right]\frac{\lambda_1} {\rho_1(a_1-1)}\cdot t^{-a_1+1}L(t).
\end{eqnarray}
\end{lemma}

\proof
First, let us rewrite (\ref{Apdx-1}) as,
\begin{eqnarray}
\gamma(s)=\frac 1 {\rho_1}-\frac 1 {\rho_1}\cdot\frac {1-\beta_1(s)} {s}\cdot\frac {s} {1-\alpha(s)}.\label{Apdx-1b}
\end{eqnarray}
In the following, we will divide the proof of Lemma~\ref{proposition-A1} into two parts, depending on whether $a_1>1$ is an integer or not.

\noindent{\it Case 1: Non-integer $a_1 >1$.} Suppose that $n<a_1<n+1$, $n\in\{1,2,\ldots\}$. Since $P\{T_{\beta_1}>t\}\sim t^{-a_1}L(t)$ and $(1-\rho_1)^{a_1+1}P\{T_{\alpha}>t\}\sim t^{-a_1}L(t)$, we know that $\beta_{1,n}<\infty$, $\beta_{1,n+1}=\infty$, $\alpha_{n}<\infty$ and $\alpha_{n+1}=\infty$.
Define $\beta_{1,n}(s)$ and $\alpha_{n}(s)$ in a manner similar to that in (\ref{phi1}). Therefore,
\begin{eqnarray}
\frac {1-\beta_1(s)} {s}&=&\beta_{1,1}+\sum_{k=2}^{n}\frac {\beta_{1,k}} {k!}(-s)^{k-1} +(-1)^{n}\frac{\beta_{1,n}(s)} {s},\label{(1-beta1-a)/s}\\
\frac {1-\alpha(s)} {s}&=&\alpha_{1}+\sum_{k=2}^{n}\frac {\alpha_{k}} {k!}(-s)^{k-1} +(-1)^{n}\frac{\alpha_{n}(s)} {s}.\label{(1-alph-a)/s}
\end{eqnarray}
By Lemma~\ref{Cohen},
\begin{equation}
(1-\rho_1)^{a_1+1}\alpha_{n}(s)\ \sim\ \beta_{n}(s)\ \sim\ \frac{\Gamma(a_1-n)\Gamma(n+1-a_1)} {\Gamma(a_1)} s^{a_1}L(1/s),\quad s\downarrow 0.\label{sim-alph-beta}
\end{equation}
Furthermore, it follows from (\ref{(1-alph-a)/s}) that,
\begin{eqnarray}
\frac {s} {1-\alpha(s)}&=&\frac {1/\alpha_1} {1+(1/\alpha_1)\sum_{k=2}^{n}\frac {\alpha_{k}} {k!}(-s)^{k-1} +(-1)^{n}(1/\alpha_1)\frac{\alpha_{n}(s)} {s}}\nonumber\\
&=&\frac 1 {\alpha_1} -\sum_{k=1}^{n-1}u_k s^k - (-1)^{n}\frac {\alpha_{n}(s)} {\alpha_1^2 s} +O(s^{n}),\label{(1-alph-a)/s-22}
\end{eqnarray}
where $u_1,u_2,\cdots,u_{n-1}$ are constants.
By (\ref{Apdx-1b}), (\ref{(1-beta1-a)/s}) and (\ref{(1-alph-a)/s-22}), we have,
\begin{eqnarray}
\gamma(s)&=&1+\sum_{k=1}^{n-1}e_k s^k +(-1)^{n} \cdot\frac {1} {\rho_1\alpha_1}\left[\frac {\beta_{1,1}} {\alpha_1}\cdot\frac{\alpha_{n}(s)} {s}-\frac{\beta_{1,n}(s)} {s}\right] +O(s^{n}),
\end{eqnarray}
where $e_1,e_2,\cdots,e_{n-1}$ are constants. Based on the above, we define $\gamma_{n-1}(s)$ in a manner similar to that in (\ref{phi1}). Applying (\ref{sim-alph-beta}), we have,
\begin{eqnarray}
\gamma_{n-1}(s)&\sim&\frac {1} {\rho_1\alpha_1}\left[\frac {\beta_{1,1}} {\alpha_1}\cdot\frac{\alpha_{n}(s)} {s}-\frac{\beta_{1,n}(s)} {s}\right]\nonumber\\
&\sim&\frac{\lambda_1} {\rho_1}\cdot\frac{1-\rho_1} {\rho_1} \left[\frac{1} {(1-\rho_1)^{a_1}}-1\right]\cdot\frac{\Gamma(a_1-n)\Gamma(n+1-a_1)} {(a_1-1)\Gamma(a_1-1)}  s^{a_1-1}L(1/s).
\end{eqnarray}
Then, making use of Lemma \ref{Cohen}, we complete the proof of Lemma~\ref{proposition-A1} for non-integer $a_1>1$.
\\
\\
{\it Case 2: Integer $a_{1}>1$.}
Suppose that $a_1=n\in\{2,3,\ldots\}$.  Since $P\{T_{\beta_1}>t\}\sim t^{-n}L(t)$ and $(1-\rho_1)^{n+1}P\{T_{\alpha}>t\}\sim t^{-n}L(t)$, we know that $\alpha_{n-1}<\infty$ and $\beta_{1,n-1}<\infty$, but, whether $\alpha_{n}$ or $\beta_{1,n}$ is finite or not remains uncertain. This uncertainty is essentially determined by whether $\int_x^{\infty}t^{-1}L(t)dt$ is convergent or not. Define $\widehat{\beta}_{1,n-1}(s)$ and $\widehat{\alpha}_{n-1}(s)$ in a way similar to that in (\ref{phi2}). Then,
\begin{eqnarray}
\frac {1-\beta_1(s)} {s}&=&\beta_{1,1}+\sum_{k=2}^{n-1}\frac {\beta_{1,k}} {k!}(-s)^{k-1} +(-s)^{n-1}\widehat{\beta}_{1,n-1}(s),\label{case2:(1-beta1-a)/s}\\
\frac {1-\alpha(s)} {s}&=&\alpha_{1}+\sum_{k=2}^{n-1}\frac {\alpha_{k}} {k!}(-s)^{k-1} +(-s)^{n-1}\widehat{\alpha}_{n-1}(s).\label{case2:(1-alph-a)/s}
\end{eqnarray}
By Lemma~\ref{Lemma 4.5new}, we obtain, for $x>0$,
\begin{eqnarray}
&&(1-\rho_1)^{n+1}\widehat{\alpha}_{n-1}(xs)-(1-\rho_1)^{n+1}\widehat{\alpha}_{n-1}(s)\nonumber\\
&&\quad \sim\ \widehat{\beta}_{1,n-1}(xs)-\widehat{\beta}_{1,n-1}(s)\ \sim\ -(\log x) L(1/s)/(n-1)!\quad \mbox{ as}\ s\downarrow 0.\label{case2:sim-alph-beta}
\end{eqnarray}
Furthermore, it follows from (\ref{case2:(1-alph-a)/s}) that,
\begin{eqnarray}
\frac {s} {1-\alpha(s)}&=&\frac {1/\alpha_1} {1+(1/\alpha_1)\sum_{k=2}^{n-1}\frac {\alpha_{k}} {k!}(-s)^{k-1} +(-s)^{n-1}(1/\alpha_1)\widehat{\alpha}_{n-1}(s)}\nonumber\\
&=&\frac 1 {\alpha_1} -\sum_{k=1}^{n-1}u'_k s^k - (-s)^{n-1}\frac {\widehat{\alpha}_{n-1}(s)} {\alpha_1^2} +O(s^{n}).\label{case2:(1-alph-a)/s-22}
\end{eqnarray}
where $u'_1,u'_2,\cdots,u'_{n-1}$ are constants.
By (\ref{Apdx-1b}), (\ref{case2:(1-beta1-a)/s}) and (\ref{case2:(1-alph-a)/s-22}), we have,
\begin{eqnarray}
\gamma(s)&=&1+\sum_{k=1}^{n-1}e'_k s^k +(-1)^{n} s^{n-1}\cdot\frac {1} {\rho_1\alpha_1}\left[\frac {\beta_{1,1}} {\alpha_1}\cdot\widehat{\alpha}_{n-1}(s)- \widehat{\beta}_{1,n-1}(s)\right] +O(s^{n}),\label{case2gamma(s)}
\end{eqnarray}
where $e'_1,e'_2,\cdots,e'_{n-1}$ are constants. Based on which, we define $\widehat{\gamma}_{n-2}(s)$ in a way similar to that in (\ref{phi2}). Then,
\begin{eqnarray}
\widehat{\gamma}_{n-2}(s)=(-1)^{n} e'_{n-1}+\frac {1} {\rho_1\alpha_1}\left[\frac {\beta_{1,1}} {\alpha_1}\cdot\widehat{\alpha}_{n-1}(s)- \widehat{\beta}_{1,n-1}(s)\right] +O(s).
\end{eqnarray}
It follows from (\ref{case2gamma(s)}) and (\ref{case2:sim-alph-beta}) that
\begin{eqnarray}
\lim_{s\downarrow 0}\frac {\widehat{\gamma}_{n-2}(xs)-\widehat{\gamma}_{n-2}(s)} {L(1/s)/(n-2)!}
&=&\frac {1} {\rho_1\alpha_1}\left[\frac {\beta_{1,1}} {\alpha_1}\cdot\lim_{s\downarrow 0}\frac {\widehat{\alpha}_{n-1}(xs)-\widehat{\alpha}_{n-1}(s)} {(n-1)L(1/s)/(n-1)!}
- \lim_{s\downarrow 0}\frac {\widehat{\beta}_{1,n-1}(xs)-\widehat{\beta}_{1,n-1}(s)} {(n-1)L(1/s)/(n-1)!} \right]\nonumber\\
&=&\frac{\lambda_1} {\rho_1}\frac{1-\rho_1} {\rho_1} \left[\frac{1} {(1-\rho_1)^{n}}-1\right]\cdot \left(-\frac 1 {n-1}\log x\right).
\end{eqnarray}
Applying Lemma \ref{Lemma 4.5new}, we complete the proof of Lemma~\ref{proposition-A1} for integer $a_1=n\in\{2,3,\ldots\}$.
\pend

\bigskip
\noindent\underline{{\sc Proof} of Proposition~\ref{lemma-Hbeta12>j}:} It follows directly from Remark~\ref{gamma-compound-possion}, Lemma~\ref{proposition-A1} and Lemma~\ref{lemma-Asmussen-poisson}. \pend

Referring to Remark~\ref{BoundAsymp}, we know from (\ref{lemma-Hbeta12>j-formula}) that $c=(1-\rho_1) \rho_1^{-1} \left[1/(1-\rho_1)^{a_1}-1\right]$. Now let us confirm that $a_1<c<a_1/(1-\rho_1)^{a_1-1}$, which is equivalent to checking that
$a_1\rho_1(1-\rho_1)^{a_1-1}+(1-\rho_1)^{a_1}<1$ and $a_1\rho_1+(1-\rho_1)^{a_1}>1$. This is true because
$a_1\rho_1(1-\rho_1)^{a_1-1}+(1-\rho_1)^{a_1}$ is decreasing in $p_1\in (0,1)$ and $a_1\rho_1+(1-\rho_1)^{a_1}$ is increasing in $p_1\in (0,1)$.

\subsubsection{Asymptotic tail probability of the r.v. $H_{\beta_2,2}$}

As we shall see in the next subsection, our main results do not require a detailed asymptotic expression for $P\{H_{\beta_2,2}>j\}$. It is enough to verify that it is $o(1)\cdot j^{-a_1+1}L(j)$ as $j\to\infty$.

Taking $z_1\to 1$ in (\ref{H{beta_2}-3}), we have,
\begin{eqnarray}
E(z_2^{H_{\beta_2,2}})=H_{\beta_2}(1,z_2)&=&\frac {p} {\rho_2}\cdot\sum_{k=1}^{\infty} k b_{\beta_2,k}\cdot D_{k}(1,z_2).\label{Hbeta2(1,z)}
\end{eqnarray}
It follows from (\ref{Hbeta2(1,z)}) that
$H_{\beta_1,2}\stackrel{\rm d}{=}D_{k,2}$, with probability $(p/\rho_2)k b_{\beta_1,k}$, for $k\ge 1$. Define the r.v. $H^U_{\beta_2,2}\stackrel{\rm d}{=}D^U_{k,2}$,   with probability $(p/\rho_2)k b_{\beta_2,k}$, for $k\ge 1$. Then, by (\ref{DDD}), we have, $H_{\beta_2,2}\le_{st}H_{\beta_2,2}^{U}$.
Note that $H_{\beta_2,2}^{U}$ has the PGF
\begin{eqnarray}
E(z_2^{H_{\beta_2,2}^{U}})&=&\frac {p} {\rho_2}\cdot\sum_{k=1}^{\infty} k b_{\beta_2,k}\cdot E(z_2^{D_{k,2}^{U}})
=\frac {p} {\rho_2}\cdot\sum_{k=1}^{\infty} k b_{\beta_1,k}\cdot (qh(z_2) +p z_2)^{k-1}. \label{Hbeta2U}
\end{eqnarray}

Let $N^*$ be a r.v. with probability distribution $P\{N^*=k\}=(p/\rho_2) k b_{\beta_2,k}$, $k\ge 1$. Therefore, (\ref{Hbeta2U}) implies
$H_{\beta_2,2}^{U}\stackrel{\rm d}{=}\sum_{k=1}^{N^*-1} Z_k$,
where $N^*$ is independent of $Z_k$, $k\ge 1$. Similar to (\ref{P{N>m}}), we can write,
\begin{equation}
P\{N^*>m\}=(p/\rho_2)\left[m\overline{b}_{\beta_2,m+1}+\sum_{k=m+1}^{\infty} \overline{b}_{\beta_2,k}\right],\label{P{N^*>m}}
\end{equation}
where $\overline{b}_{\beta_2,k}=\sum_{n=k}^{\infty} b_{\beta_2,n}$.
By the definition of $b_{\beta_2,n}$ in Subsection \ref{subsection-M2} and applying Lemma~\ref{lemma-Asmussen-poisson} and Lemma~\ref{lemma-Grandell-poisson}, we know that $\overline{b}_{\beta_2,k}=P\{N_{\lambda}(T_{\beta_2})>k-1\}=O(1)\cdot k^{-a_2}L(k)$ as $k\to\infty$. Furthermore, by (\ref{P{N^*>m}}) and applying Proposition~1.5.10 in \cite{Bingham:1989}, we have,
\begin{eqnarray}
P\{N^*>m\}&=& O(1)\cdot  m^{-a_2+1}L(m)\quad{as}\ m\to\infty.
\end{eqnarray}

As pointed out in Subsection \ref{subsectionH12},
$P\{Z_k>j\}= Const\cdot j^{-a_1} L(j)$ as $j\to\infty$. By Lemma~\ref{Lemma 4.5}, we know
$P\{H_{\beta_2,2}^{U}>j\}= O(1)\cdot\max\left(j^{-a_2+1}L(j),j^{-a_1}L(j)\right)$ as $j\to\infty$.
Since $P\{H_{\beta_2,2}>j\}\le P\{H_{\beta_2,2}^{U}>j\}$ and $a_2>a_1$, we have,
\begin{eqnarray}
P\{H_{\beta_2,2}>j\}&=&O(1)\cdot\max\left(j^{-a_2+1}L(j),j^{-a_1}L(j)\right)=o(1)\cdot j^{-a_1+1}L(j).\label{Hbeta22-tail}
\end{eqnarray}

\subsubsection{Asymptotic tail probabilities of the r.v.s $R_{ik},\ i,k=1,2$}

We first provide tail asymptotic probabilities for the r.v.s $M_{ik}$, $i,k=1,2$.
By (\ref{M1k>j}) and applying Lemma~\ref{Embrechts-compound-geo}, together with (\ref{Sbeta-i1-tail}), we have,
\begin{eqnarray}
P\{M_{11}>j\}&\sim&\frac {\rho_1} {1-\rho_1} P \{H_{\beta_1,1}>j \}=\frac {\rho_1} {1-\rho_1} P \{S_{\beta_1,1}>j \}
\quad\mbox{(refer to (\ref{Sbeta11-tail}))}, \label{M11-tail-a}\\
P\{M_{12}>j\}&\sim&\frac {\rho_1} {1-\rho_1} P \{H_{\beta_1,2}>j \}
\quad\mbox{(refer to (\ref{lemma-Hbeta12>j-formula}))}. \label{M12-tail-a}
\end{eqnarray}

Immediately, from (\ref{M21>j}) and (\ref{Sbeta-i1-tail}),
\begin{eqnarray}
P\{M_{21}>j\}&=&\vartheta P \{H_{\beta_2,1} >j \} =\vartheta P \{S_{\beta_2,1} >j \}
\quad\mbox{(refer to (\ref{Sbeta21-tail}))}.\label{M21-tail-a}
\end{eqnarray}
By (\ref{M22>j}) and applying Lemma~\ref{Lemma 3.4new}, together with (\ref{Hbeta22-tail}),
\begin{eqnarray}
P\{M_{22}>j\}&=& \vartheta P\{H_{\beta_2,2}+K_a+K_c>j\}\sim  \vartheta P \{K_a+K_c>j \}
\quad\mbox{(refer to (\ref{Ka+Kc-asym}))}.\label{M22-tail-a}
\end{eqnarray}

Now we are in the position to present the tail asymptotic probabilities for the r.v.s $R_{ik},\ i,k=1,2$. Recall (\ref{R11-decom}) and (\ref{R12-decom}).
By (\ref{M21-tail-a}) and (\ref{P{R_0>j}sim}), $M_{21}$ and $R_0$ have tail probabilities lighter than $j^{-a_1+1}L(j)$, and by
(\ref{M11-tail-a}), (\ref{M22-tail-a}), (\ref{Sbeta11-tail}) and (\ref{Sbeta12-tail}),
$M_{11}$, $M_{22}$, $S_{\beta_1,1}$ and $S_{\beta_1,2}$ have regularly varying tails with index $-a_1+1$.
Applying Lemma~\ref{Lemma 3.4new}, we obtain,
\begin{eqnarray}
P\{R_{11}>j\}&=&P\{M_{21}+M_{11}+S_{\beta_1,1}>j\}\nonumber\\
&\sim& P\{M_{11}+S_{\beta_1,1}>j\}\sim\frac {1} {1-\rho_1} P \{S_{\beta_1,1}>j \}
\quad\mbox{(refer to (\ref{Sbeta11-tail}))},\\
P\{R_{12}>j\}&=&P\{M_{22}+M_{12}+S_{\beta_1,2}+R_0>j\}\nonumber\\
&\sim& P\{M_{22}+M_{12}+S_{\beta_1,2}>j\}\quad\mbox{(refer to (\ref{M22-tail-a}), (\ref{M12-tail-a}) and (\ref{Sbeta12-tail}))}. \label{4.48}
\end{eqnarray}

By a similar argument, it follows from (\ref{R21-decom})--(\ref{R22-decom}) that,
\begin{eqnarray}
P\{R_{21}>j\}&=&P\{S_{\beta_2,1}>j\}\quad\mbox{(refer to (\ref{Sbeta21-tail}))},\\
P\{R_{22}>j\}&=&P\{S_{\beta_2,2}+ K_a+ K_c +R_0>j\}\sim P\{K_a+ K_c>j\}\quad\mbox{(refer to (\ref{Ka+Kc-asym}))},
\end{eqnarray}
where we have used the fact, by (\ref{Sbeta22-tail}), that $S_{\beta_2,2}$ has a tail probability lighter than $j^{-a_1+1}L(j)$.

Recall the definition of $R_{i,k}$, $i,k=1,2$ in Section 2. We know that $P\{R_{que}>j|I_{ser}=i\}=P\{R_{i1}>j\}$ and $P\{R_{orb}>j|I_{ser}=i\}=P\{R_{i2}>j\}$, $i=1,2$. The above discussion is summarized in the following theorem.
\begin{theorem} \label{The:4.2}
As $j\to\infty$,
\begin{eqnarray}
P\{R_{que}>j|I_{ser}=1\}&\sim& \frac {\lambda_1^{a_1}} {\rho_1 (1-\rho_1) (a_1-1) } \cdot j^{-a_1+1} L(j),\label{R11-tail-final}\\
P\{R_{orb}>j|I_{ser}=1\}&\sim& \left[\frac {\rho_2}{1-\rho}+ \frac{1} {\rho_1}\right]\cdot\frac{\lambda_1 \lambda_2^{a_1-1}} {(a_1-1)(1-\rho_1)^{a_1}}\cdot j^{-a_1+1}L(j),\label{R12-tail-final}\\
P\{R_{que}>j|I_{ser}=2\}&\sim& \left\{\begin{array}{ll}
\frac {\lambda_2\lambda_1^{a_2-1}} {\rho_2 (a_2-1)}\cdot j^{-a_2+1} L(j), &\mbox{ if }r=0,\\
\frac {\lambda_2\lambda_1 (\lambda_1+r)^{a_2-1}} {\rho_2 r} \cdot \left(\frac {\lambda_1} {\lambda_1+r}\right)^j j^{-a_2} L(j), &\mbox{ if }r>0,
\end{array}
\right.\label{R21-tail-final}\\
P\{R_{orb}>j|I_{ser}=2\}&\sim& \frac {\lambda_1\lambda_2^{a_1-1}}{(1-\rho) (a_1-1) (1-\rho_1)^{a_1-1}} \cdot j^{-a_1+1} L(j).\label{R22-tail-final}
\end{eqnarray}
\end{theorem}

To conclude the paper, we would like to provide intuition on the results in Theorem~\ref{The:4.2}. First, let us recall a well-known result for the standard $M/G/1$ queue: if the service time is regularly varying with index $-a_1$, then the stationary queue length is also regularly varying, but with index $-a_1+1$. Such a conclusion can be made through a distributional Little's law (see, e.g., \cite{Asmussen-Klupperlberg-Sigman:1999}). For the model studied in this paper, the condition $I_{ser}=1$ means that the server is serving a Type-1 customer. Under this condition, both types of customers have to wait, customers of Type-1 in the queue and customers of Type-2 in the orbit. Therefore, both $R_{que}|I_{ser}=1$ and $R_{orb}|I_{ser}=1$ have the asymptotic tail in the form of $Const\cdot j^{-a_1+1}L(j)$ (given in (\ref{R11-tail-final}) and (\ref{R12-tail-final})), due to the regularly varying assumption for the service time of Type-1 customers in Assumption~A1.
On the other hand, the condition $I_{ser}=2$ means that the server is serving a Type-2 (lower priority) customer, which implies that no Type-1 customers were waiting in the queue at the beginning of service of this Type-2 customer. In other words, $I_{ser}=2$ implies that all Type-1 customers in the queue must be those who arrived after the beginning of the service time of this Type-2 customer. Therefore,
$R_{que}|I_{ser}=2$ has an asymptotic tail in the form given in (\ref{R21-tail-final}), determined by the service time assumption (in Assumption~A2) of Type-2 customers.
However, $R_{orb}|I_{ser}=2$ still has an asymptotic tail in the form of $Const\cdot j^{-a_1+1}L(j)$ (by (\ref{R22-tail-final})) (determined by the assumption on the Type-1 customer's service time), since the customers arrived to the orbit could be those arrived during the service times of Type-2 customers, and/or Type-1 customers who were served before the current Type-2 customer in service, due to the priority discipline, and the tail of the service time for Type-1 customers is heavier than that for Type-2 customers.

\section*{Acknowledgments}
This work was supported in part by the National Natural Science Foundation of China (Grant No. 71571002),
the Natural Science Foundation of the Anhui Higher Education Institutions of China (No. KJ2017A340),
the Research Project of Anhui Jianzhu University, and a Discovery Grant from the Natural Sciences and Engineering Research Council of Canada (NSERC).

\appendix

\section{Appendix}

\subsection{Definitions and useful results from the literature}

\begin{definition}[for example, see Bingham, Goldie and Teugels~\cite{Bingham:1989}]
\label{Definition 3.1}
A measurable function $U:(0,\infty)\to (0,\infty)$ is regularly varying at $\infty$ with index $\sigma\in(-\infty,\infty)$ (written $U\in \mathcal{R}_{\sigma}$) iff $\lim_{t\to\infty}U(xt)/U(t)=x^{\sigma}$ for all $x>0$. If $\sigma=0$ we call $U$ slowly varying, i.e., $\lim_{t\to\infty}U(xt)/U(t)=1$ for all $x>0$.
\end{definition}
\begin{definition}[for example, see Foss, Korshunov and Zachary~\cite{Foss2011}]
\label{Definition 3.2}
A distribution $F$ on $(0,\infty)$ belongs to the class of subexponential distribution (written $F\in \mathcal S$) if $\lim_{t\to\infty}\overline{F^{*2}}(t)/\overline{F}(t)=2$, where $\overline{F}=1-F$ and $F^{*2}$ denotes the second convolution of $F$.
\end{definition}

\begin{lemma}[de Meyer and Teugels~\cite{Meyer-Teugels:1980}]
\label{Lemma 2.1}
Under Assumption A1,
\begin{equation}\label{Meyer-Teugels-T{alpha}}
P\{T_{\alpha}>t\}\sim \frac 1 {(1-\rho_1)^{a_1+1}}\cdot t^{-a_1} L(t)\quad \mbox{as }t\to\infty.
\end{equation}
\end{lemma}
The result (\ref{Meyer-Teugels-T{alpha}}) is straightforward due to the main theorem in \cite{Meyer-Teugels:1980}.

\begin{lemma}[pp.580--581 in \cite{Embrechts1997}]
\label{Embrechts-compound-geo}
Let $N$ be a r.v. with $P\{N=k\}=(1-\sigma)\sigma^{k-1}$, $0<\sigma<1$, $k\ge 1$, and
$\{Y_k\}_{k=1}^{\infty}$ be a sequence of non-negative, i.i.d. r.v.s having a common subexponential distribution $F$.
Define $S_n=\sum_{k=1}^n Y_k$. Then $P\{S_N > t\} \sim (1-F(t))/(1-\sigma)$ as $t\to \infty$.
\end{lemma}

\begin{lemma}[Proposition~3.1 in  \cite{Asmussen-Klupperlberg-Sigman:1999}]
\label{lemma-Asmussen-poisson} Let $N_{\lambda}(t)$ be a Poison process with rate $\lambda$ and let $T$ be a positive r.v. with distribution $F$, which is independent of $N_{\lambda}(t)$.
If $\bar{F}(t)=P\{T>t\}$ is heavier than $e^{-\sqrt{t}}$ as $t\to\infty$, then $P(N_{\lambda}(T)>j)\sim P\{T>j/\lambda\}$ as $j\to\infty$.
\end{lemma}

Lemma~\ref{lemma-Asmussen-poisson} holds for any distribution $F$ with a regularly varying tail because it is heavier than $e^{-\sqrt{t}}$ as $t\to\infty$.

\begin{lemma}[p.181 in \cite{Grandell:1997}]
\label{lemma-Grandell-poisson}
Let $N_{\lambda}(t)$ be a Poison process with rate $\lambda$ and let $T$ be a positive r.v. with distribution $F$, which is independent of $N_{\lambda}(t)$. If $\bar{F}(t)\stackrel{\rm def}{=}P\{T>t\}\sim e^{-w t} t^{-h}L(t)$ as $t\to\infty$ for $w> 0$ and $-\infty<h<\infty$, then
\[
P(N_{\lambda}(T)>j)\sim \lambda (\lambda+w)^{h-1}\left(\frac {\lambda} {\lambda+w}\right)^j j^{-h}L(j),\quad j\to \infty.
\]
\end{lemma}
\begin{lemma}[p.48 in \cite{Foss2011}]
\label{Lemma 3.4new}
Let $F$, $F_1$ and $F_2$ be distribution functions.
Suppose that $F\in\mathcal S$.
If $\bar{F}_i(t)/\bar{F}(t)\to c_i$ as $t\to\infty$ for some $c_i\ge 0, \; i=1,2$, then
$\overline{F_1*F}_2(t)/\bar{F}(t)\to c_1+c_2$ as $t\to\infty$,
where the symbol $\bar{F}\stackrel{\rm def}{=}1-F$ and  ``$F_1*F_2$" stands for the convolution of $F_1$ and $F_2$.
\end{lemma}
\begin{lemma}[pp.162--163 in \cite{Grandell:1997}]
\label{Lemma 4.5}
Let $N$ be a discrete non-negative integer-valued r.v. with mean value $\mu_N$,
and $\{Y_k\}_{k=1}^{\infty}$ be a sequence of non-negative i.i.d. r.v.s with mean value $\mu_Y$.
Define $S_0\equiv 0$ and $S_n=\sum_{k=1}^n Y_k$.
If $P\{Y_k>x\}\sim c_Y x^{-h} L(x) $ as $x\to\infty$ and $P\{N>m\}\sim c_N m^{-h}L(m)$ as $m\to\infty$, where $h> 1$,  $c_Y\ge 0$ and $c_N\ge 0$, then
$P\{S_N > x\}\sim (c_N\mu_Y^h + \mu_N c_Y) x^{-h}L(x)$ as $x\to \infty.$
\end{lemma}

\begin{remark}
It is a convention that in Lemma~\ref{Lemma 4.5}, $c_Y=0$ and $c_N=0$ means that
$\lim_{x\to\infty}P\{Y_k>x\}/(x^{-h} L(x))=0$ and $\lim_{m\to\infty}P\{N>m\}/(m^{-h} L(m))=0$, respectively.
\end{remark}

The following two criteria are from Feller (1991)~\cite{Feller1971} (see p.441), which are often used to verify that a function is completely monotone.

\vskip 0.5cm

\noindent{\bf Criterion A.1} If $\vartheta_1(\cdot)$ and $\vartheta_2(\cdot)$ are completely monotone so is their product $\vartheta_1(\cdot)\vartheta_2(\cdot)$.

\vskip 0.5cm

\noindent{\bf Criterion A.2} If $\vartheta_3(\cdot)$ is completely monotone and $\vartheta_4(\cdot)$ a positive function with a completely monotone
derivative $\vartheta'_4(\cdot)$ then $\vartheta_3(\vartheta_4(\cdot))$ is completely monotone.

\vskip 0.5cm

To prove Lemma~\ref{proposition-A1}, let us list some notations and results, which will be used.
Let $F(x)$ be any distribution on $[0,\infty)$ with the LST $\phi(s)$.
We denote the $n$th moment of $F(x)$ by $\phi_n$, $n\ge 0$.
It is well known that $\phi_n<\infty$ iff
\begin{equation}\label{phi00}
\phi(s)=\sum_{k=0}^{n}\frac{\phi_k}{k!}(-s)^k + o(s^n),\quad n\ge 0.
\end{equation}
Based on (\ref{phi00}), we introduce the notation $\phi_n(s)$ and $\widehat{\phi}_n(s)$, defined by
\begin{eqnarray}
\phi_n(s)&\stackrel{\rm def}{=}&(-1)^{n+1}\left\{\phi(s)-\sum_{k=0}^{n}\frac{\phi_k}{k!}(-s)^k\right\},\quad n\ge 0,\label{phi1}\\
\widehat{\phi}_n(s)&\stackrel{\rm def}{=}&\phi_n(s)/s^{n+1},\quad n\ge 0.\label{phi2}
\end{eqnarray}
\begin{lemma}[pp.333--334 in \cite{Bingham:1989}]
\label{Cohen}
Assume that $n<d<n+1$, $n\in\{0,1,2,\ldots\}$, then the following two statements are equivalent:
\begin{eqnarray}
&&1-F(t) \sim  t^{-d} L(t),\quad t\to\infty;\\
&&\phi_n(s) \sim \frac{\Gamma(d-n)\Gamma(n+1-d)} {\Gamma(d)} s^{d}L(1/s),\quad s\downarrow 0.
\end{eqnarray}
\end{lemma}

\begin{lemma}[Lemma 3.3 in \cite{Liu-Zhao:2017b}]
\label{Lemma 4.5new}
Assume that $n\in\{1,2,\ldots\}$, then the following two statements are equivalent:
\begin{eqnarray}
&&1-F(t)\sim t^{-n}L(t),\quad t\to\infty; \label{deHaan3}\\
&&\lim_{s\downarrow 0}\frac {\widehat{\phi}_{n-1}(xs)-\widehat{\phi}_{n-1}(s)}{L(1/s)/(n-1)!}=-\log x, \quad\mbox{ for all }x>0.\label{deHaan4}
\end{eqnarray}
\end{lemma}

In \cite{Liu-Zhao:2017b}, Lemma~\ref{Lemma 4.5new} is proved by applying Karamata's theorem in \cite{Bingham:1989}, p.27, the monotone density theorem in \cite{Bingham:1989}, p.39 and Theorem~3.9.1 in \cite{Bingham:1989}, pp.172--173.

\end{document}